\magnification=\magstep1
\hsize=16.5 true cm 
\vsize=23.6 true cm
\font\bff=cmbx10 scaled \magstep1
\font\bfff=cmbx10 scaled \magstep2

\font\bffgg=cmbx10 scaled \magstep4
\font\smc=cmcsc10 
\parindent0cm
\def\cl{\centerline}           %
\def\bp{\bigskip}              %
\def\mp{\medskip}              %
\def\sp{\smallskip}            %
\def\bc{{\bf c}}               %
           %
\def\Bbb#1{\hbox{\boldmas #1}} %
\def\Q{\Bbb Q}                 %
\def\R{\Bbb R}                 %
\def\N{\Bbb N}                 %
\def\Z{\Bbb Z}                 %

\expandafter\edef\csname amssym.def
\endcsname{%
       \catcode`\noexpand\@=\the\catcode`\@\space}

\catcode`\@=11

\def\undefine#1{\let#1\undefined}
\def\newsymbol#1#2#3#4#5{\let\next@\relax
 \ifnum#2=\@ne\let\next@\msafam@\else
 \ifnum#2=\tw@\let\next@\msbfam@\fi\fi
 \mathchardef#1="#3\next@#4#5}
\def\mathhexbox@#1#2#3{\relax
 \ifmmode\mathpalette{}{\m@th\mathchar"#1#2#3}%
 \else\leavevmode\hbox{$\m@th\mathchar"#1#2#3$}\fi}
\def\hexnumber@#1{\ifcase#1 0\or 1\or 2\or 3\or 4\or 5\or 6\or 7\or 8\or
 9\or A\or B\or C\or D\or E\or F\fi}

\font\tenmsa=msam10
\font\sevenmsa=msam7
\font\fivemsa=msam5
\newfam\msafam
\textfont\msafam=\tenmsa
\scriptfont\msafam=\sevenmsa
\scriptscriptfont\msafam=\fivemsa
\edef\msafam@{\hexnumber@\msafam}
\mathchardef\dabar@"0\msafam@39
\def\dashrightarrow{\mathrel{\dabar@\dabar@\mathchar"0\msafam@4B}}
\def\dashleftarrow{\mathrel{\mathchar"0\msafam@4C\dabar@\dabar@}}

\def\ulcorner{\delimiter"4\msafam@70\msafam@70 }
\def\urcorner{\delimiter"5\msafam@71\msafam@71 }
\def\llcorner{\delimiter"4\msafam@78\msafam@78 }
\def\lrcorner{\delimiter"5\msafam@79\msafam@79 }
\def\yen{{\mathhexbox@\msafam@55}}
\def\checkmark{{\mathhexbox@\msafam@58}}
\def\circledR{{\mathhexbox@\msafam@72}}
\def\maltese{{\mathhexbox@\msafam@7A}}

\font\tenmsb=msbm10
\font\sevenmsb=msbm7
\font\fivemsb=msbm5
\newfam\msbfam
\textfont\msbfam=\tenmsb
\scriptfont\msbfam=\sevenmsb
\scriptscriptfont\msbfam=\fivemsb
\edef\msbfam@{\hexnumber@\msbfam}
\def\Bbb#1{{\fam\msbfam\relax#1}}
\def\widehat#1{\setbox\z@\hbox{$\m@th#1$}%
 \ifdim\wd\z@>\tw@ em\mathaccent"0\msbfam@5B{#1}%
 \else\mathaccent"0362{#1}\fi}

\def\widetilde#1{\setbox\z@\hbox{$\m@th#1$}%
 \ifdim\wd\z@>\tw@ em\mathaccent"0\msbfam@5D{#1}%
 \else\mathaccent"0365{#1}\fi}
\font\teneufm=eufm10
\font\seveneufm=eufm7
\font\fiveeufm=eufm5
\newfam\eufmfam
\textfont\eufmfam=\teneufm
\scriptfont\eufmfam=\seveneufm
\scriptscriptfont\eufmfam=\fiveeufm

\newsymbol\risingdotseq 133A
\newsymbol\fallingdotseq 133B
\newsymbol\complement 107B
\newsymbol\nmid 232D
\newsymbol\rtimes 226F
\newsymbol\thicksim 2373

\font\eightmsb=msbm8   \font\sixmsb=msbm6   \font\fivemsb=msbm5
\font\eighteufm=eufm8  \font\sixeufm=eufm6  \font\fiveeufm=eufm5
\font\eightrm=cmr8     \font\sixrm=cmr6     \font\fiverm=cmr5
\font\eightbf=cmbx8    \font\sixbf=cmbx6    
      \font\eighti=cmmi8   \font\sixi=cmmi6
\font\ninesy=cmsy9     \font\eightsy=cmsy8  \font\sixsy=cmsy6
     \font\eightit=cmti8  
     \font\eightsl=cmsl8  
     \font\eighttt=cmtt8

\font\eightsmc=cmcsc8
\newskip\ttglue
\newfam\smcfam
\def\eightpoint{\def\rm{\fam0\eightrm}%
  \textfont0=\eightrm \scriptfont0=\sixrm \scriptscriptfont0=\fiverm
  \textfont1=\eighti \scriptfont1=\sixi \scriptscriptfont1=\fivei
  \textfont2=\eightsy \scriptfont2=\sixsy \scriptscriptfont2=\fivesy
  \textfont3=\tenex \scriptfont3=\tenex \scriptscriptfont3=\tenex
  \def\smc{\fam\smcfam\eightsmc}
  \textfont\smcfam=\eightsmc          
\textfont\eufmfam=\eighteufm              \scriptfont\eufmfam=\sixeufm
     \scriptscriptfont\eufmfam=\fiveeufm
\textfont\msbfam=\eightmsb            \scriptfont\msbfam=\sixmsb
     \scriptscriptfont\msbfam=\fivemsb
\def\it{\fam\itfam\eightit}%
  \textfont\itfam=\eightit
  \def\sl{\fam\slfam\eightsl}%
  \textfont\slfam=\eightsl
  \def\bf{\fam\bffam\eightbf}%
  \textfont\bffam=\eightbf \scriptfont\bffam=\sixbf
   \scriptscriptfont\bffam=\fivebf
  \def\tt{\fam\ttfam\eighttt}%
  \textfont\ttfam=\eighttt
  \tt \ttglue=.5em plus.25em minus.15em
  \normalbaselineskip=9pt
  \def\MF{{\manual opqr}\-{\manual stuq}}%
  \let\big=\eightbig
  \setbox\strutbox=\hbox{\vrule height7pt depth2pt width\z@}%
  \normalbaselines\rm}
\def\eightbig#1{{\hbox{$\textfont0=\ninerm\textfont2=\ninesy
  \left#1\vbox to6.5pt{}\right.\n@space$}}}

\catcode`@=13 

\cl{\bffgg Counting continua}
\bp
\cl{\bfff Gerald Kuba}
\bp
\vbox{\eightpoint
{\bf Abstract.}  For infinite cardinals $\,\kappa,\lambda\,$
let $\,{\cal C}(\kappa,\lambda)\,$ denote the class
of all {\it compact Hausdorff spaces} of {\it weight} $\,\kappa\,$ and 
{\it size} $\,\lambda\,$.
So $\,{\cal C}(\kappa,\lambda)=\emptyset\,$
if $\,\kappa>\lambda\,$ or $\,\lambda>2^\kappa\,$. 
If $\,{\cal F}\,$ is a class of pairwise non-homeomorphic spaces 
in $\,{\cal C}(\kappa,\lambda)\,$ 
then $\,{\cal F}\,$ is a set of size not greater
than $\,2^\kappa\,$.  
For every infinite cardinal $\,\kappa\,$ we construct 
$\,2^\kappa\,$ pairwise non-embeddable {\it pathwise connected}
spaces in $\,{\cal C}(\kappa,\lambda)\,$ 
for $\,\lambda=\max\{2^{\aleph_0},\kappa\}\,$ and for
$\,\lambda=\exp\log(\kappa^+)\,$. 
(If $\,\kappa\,$ is a strong limit then $\,\exp\log(\kappa^+)=2^\kappa\,$.)
Additionally, for all 
infinite cardinals $\,\kappa,\mu\,$ with $\,\mu\leq\kappa\,$
we construct $\,2^\kappa\,$ pairwise non-embeddable {\it connected}
spaces in $\,{\cal C}(\kappa,\kappa^\mu)\,$. 
Furthermore, for $\,\kappa=\lambda=2^\theta\,$ with arbitrary $\,\theta\,$
and for certain other pairs $\,\kappa,\lambda\,$
we construct $\,2^{\kappa}\,$ pairwise non-embeddable connected,
{\it linearly ordered} spaces $\,X\in{\cal C}(\kappa,\lambda)\,$ such that 
$\,Y\in{\cal C}(\kappa,\lambda)\,$
whenever $\,Y\,$ is an infinite compact and connected subspace of $\,X\,$.
On the other hand we prove that there is no space $\,X\,$
with this property either if $\,\kappa=\lambda\,$ and 
$\,{\rm cf}\,\kappa=\aleph_0\,$ or if $\,\lambda\,$ is a strong limit
with $\,{\rm cf}\,\lambda=\aleph_0\,$. 
\sp   
}
\bp
{\bff 1. Introduction}
\mp
Write $\,|S|\,$ for the cardinality ({\it size}) of a set $\,S\,$.
Cardinals are the initial ordinals and
$\,\kappa,\lambda,\mu\,$ are used throughout to stand for 
infinite cardinals and $\,\bc=|\R|=2^{\aleph_0}\,$ is the cardinality of the 
continuum and $\,\kappa^+\,$ is the least cardinal
greater than $\,\kappa\,$, whence $\,\aleph_1=\aleph_0^+\leq\bc\,$.
For $\,\lambda\leq\kappa\,$ 
let $\,\nu(\lambda,\kappa)\,$ denote the smallest $\,\mu\,$
such that $\,\lambda^\mu>\kappa\,$.
In particular, if $\,\lambda\leq \bc\;$
then $\;\nu(\lambda,\kappa)=\log(\kappa^+)=
\min\{\,\mu\;|\;2^\mu>\kappa\,\}\;$
and hence $\;\lambda^{\nu(\lambda,\kappa)}=2^{\nu(\bc,\kappa)}\,$.
Of course, if $\,\lambda\leq\kappa<\bc\,$ then 
$\;\nu(\lambda,\kappa)=\aleph_0\;$ and
$\;\lambda^{\nu(\lambda,\kappa)}=\bc\,$.
It is plain that if $\;\lambda_1\leq \lambda_2\leq\kappa\;$ then   
\sp
\cl{$\,\nu(\lambda_2,\kappa)\leq\nu(\lambda_1,\kappa)\leq\kappa<\kappa^+\leq
\lambda_2^{\nu(\lambda_2,\kappa)}\leq \lambda_1^{\nu(\lambda_1,\kappa)}
=\kappa^{\nu(\lambda_1,\kappa)}\leq 2^\kappa\,$.} 
\sp
Note that (by [3] 5.20)
if $\,\kappa\,$ is a strong limit 
(i.e.~$\,2^\theta<\kappa\,$ for every cardinal $\,\theta<\kappa\,$)
then $\,\nu(\kappa,\kappa)\,$ is the cofinality of $\,\kappa\,$ and hence
$\;\lambda^{\nu(\lambda,\kappa)}=2^\kappa\;$
for every $\,\lambda\leq\kappa\,$.
On the other hand, for certain $\,\kappa>\bc\,$ one cannot rule out
$\;\kappa^{\nu(\kappa,\kappa)}<\bc^{\nu(\bc,\kappa)}<2^\kappa\;$
even with extremely large gaps between the three cardinalities.
Furthermore, for every regular $\,\kappa\,$  
the equation 
$\,|\{\,\lambda\;|\;\kappa<\lambda<\kappa^{\nu(\kappa,\kappa)}\,\}|=
\kappa^{\nu(\kappa,\kappa)}=2^\kappa\,$ 
is consistent with ZFC set theory.
And it is consistent with ZFC set theory that 
$\,|\{\,\lambda\;|\;\kappa<\lambda<2^\kappa\}|=2^\kappa\,$ 
holds for every regular $\,\kappa\,$, particularly for 
$\,\kappa=\aleph_0\,$.
(For explanations of these statements see Section 11.)
\mp   
Referring to [2], [4]
we use the standard notations $\,w(X),d(X),\chi(X)\,$
for {\it weight, density, character} of a space $\,X\,$.
Note that $\,\max\{\chi(X),d(X)\}\leq w(X)\leq |X|\leq 2^{\chi(X)}\,$
for every compact Hausdorff space $\,X\,$.
The following estimate is essential.
\sp
(1.1)$\;$ {\it If $\,{\cal F}\,$ 
is any family of mutually non-homeomorphic
compact Hausdorff spaces of weight $\,\kappa\,$ 
then $\;|{\cal F}|\leq 2^\kappa\,$.}
\sp
(1.1) is true because every normal space of weight $\,\kappa\,$ is
embeddable into the Hilbert cube $\,[0,1]^\kappa\,$ and
the compact Hausdorff space $\,[0,1]^\kappa\,$ 
has precisely $\,2^\kappa\,$ closed subspaces.
In (1.1) the upper bound $\,2^\kappa\,$ 
can be achieved for every $\,\kappa\,$. 
Moreover, $\,|{\cal F}|=2^\kappa\,$ can always be achieved 
with the restriction that every space in $\,{\cal F}\,$ is {\it connected}, 
see [9]. 
Alternatively, 
$\,|{\cal F}|=2^\kappa\,$ can be achieved 
for every $\,\kappa>\aleph_0\,$ such that every $\,X\in{\cal F}\,$
is {\it scattered} (and hence $|X|=w(X)=\kappa$), see [6].
(This result is unprovable for $\,\kappa=\aleph_0\,$ in view of [10].)
\mp
By analogy with a standard definition in [11],
we call spaces $\;X_i\;(i\in I)\;$ {\it incomparable} if and only if 
$\,X_i\,$ is not homeomorphic to a subspace of $\,X_j\,$
whenever $\,i,j\in I\,$ and $\,i\not=j\,$.
(In other words, spaces are incomparable if they are pairwise
non-embeddable.)
A fortiori, incomparable spaces are mutually non-homeomorphic.
If $\,X,Y\,$ are normal spaces of weight $\,\kappa\,$
and $\,X\,$ contains a copy of the
Hilbert cube $\,[0,1]^\kappa\,$ then $\,X,Y\,$ are {\it not incomparable}.
Fortunately, we can overcome this obstacle and prove the following theorem.
\eject
\mp\sp
{\bf Theorem 1.} {\it If $\,\kappa\,$ is an infinite cardinal 
and $\,\theta=\max\{\kappa,\bf c\}\,$ or
$\,\theta=\lambda^{\nu(\lambda,\kappa)}\,$
for $\,\lambda\leq\kappa\,$ 
then there exist $\,2^\kappa\,$ incomparable pathwise connected,
compact Hausdorff spaces
of weight $\,\kappa\,$ and size $\,\theta\,$.
For $\,\theta=\bc\,$ it can be accomplished
that all spaces are first countable.
Furthermore, if $\,\mu\leq\kappa\,$ then there exist 
$\,2^\kappa\,$ incomparable connected,
compact Hausdorff spaces
of weight $\,\kappa\,$ and size $\,\kappa^\mu\,$.}
\mp
Concerning the last statement in Theorem 1 
note that $\,\kappa^\kappa=2^\kappa\,$.
In the important case $\,\kappa=\aleph_0\,$ 
the statement in Theorem 1 can be specified in the following noteworthy way.
\mp
{\bf Theorem 2.} {\it There exist $\,\bc\,$ incomparable,
pathwise connected, compact subspaces of $\,\R^3\,$.}
\mp                         
In the following, a {\it continuum} 
is any compact and connected Hausdorff space 
with more than one point. Thus $\,|X|\geq \bc\,$ for any continuum $\,X\,$.
A {\it linearly ordered space} $\,X\,$ is a space
whose topology is the order topology of some
linear ordering  of $\,X\,$. (Just as metric spaces, linearly ordered spaces
are hereditarily normal.) 
In order to prove the first statement of Theorem 1 for $\,\kappa>\aleph_0\,$
by constructing appropriate {\it pathwise connected} spaces 
we will apply the following theorem about 
{\it totally pathwise disconnected} spaces. 
\mp
{\bf Theorem 3.} {\it If $\,\kappa>\aleph_0\,$ 
and $\,\lambda\leq \kappa\,$ then
there exist families 
$\,{\cal F}_\kappa,{\cal G}_\lambda\,$ 
of incomparable,
totally pathwise disconnected linearly ordered continua such that 
$\,|{\cal F}_\kappa|=|{\cal G}_{\lambda}|=2^\kappa\,$
and $\,w(X)=\kappa\,$ for every space $\,X\,$ in 
$\,{\cal F}_\kappa\cup{\cal G}_\lambda\,$
and $\,|X|=\max\{\kappa,\bc\}\,$ for every $\,X\in{\cal F}_\kappa\,$
and $\,|X|=\lambda^{\nu(\lambda,\kappa)}\,$ 
(and hence $\,w(X)<|X|\,$) for every 
$\,X\in{\cal G}_\lambda\,$.}
\mp
In Theorem 3 the case $\,\kappa=\aleph_0\,$ is excluded
since any second countable linearly ordered continuum is 
homeomorphic to the unit interval [0,1] and hence
pathwise connected, see (6.1) below.
(Conversely, it is trivial that every pathwise connected
linearly ordered continuum is homeomorphic to [0,1].)
For linearly ordered continua of size $\,\bc\,$ 
Theorem~3 can be sharpened in the following way.
(Note that $\,|X|=\bc\,$ if $\,X\,$ is a first countable 
compact Hausdorff space and either $\,w(X)>\aleph_0\,$
or $\,X\,$ is dense in itself.)
\mp
{\bf Theorem 4.} {\it For every uncountable $\,\kappa\leq \bc\,$ there exist 
$\,2^\kappa\,$ incomparable, totally pathwise disconnected, 
first countable linearly ordered continua of weight $\,\kappa\,$.}
\mp
In comparison (see [6] Theorem 5), for every $\,\kappa\leq \bc\,$ 
there exist $\,2^\kappa\,$ mutually non-homeo\-morphic separable,
dense-in-itself, compact linearly ordered spaces of
weight $\,\kappa\,$. (Note that any {\it separable}
linearly ordered space is first countable.)
The challenge in proving Theorem~4 is the case $\,\kappa=\aleph_1\,$
in view of [1] and the following theorem.
\mp
{\bf Theorem 5.} {\it It is consistent with {\rm ZFC}
that there exist precisely three spaces $\,X\,$ up
to homeomorphism such that $\,X\,$ is a first countable 
linearly ordered continuum of weight $\,\aleph_1\,$
where the union of all singleton path components is a 
separable subspace of $\,X\,$.}
\mp
Following [3], if $\,\theta\geq 2\,$ is a cardinal
then $\;\theta^{<\mu}\;$ denotes the supremum of 
all cardinals $\,\theta^\eta\,$ where $\,\eta\,$
runs through the finite and infinite cardinals smaller than $\,\mu\,$.
So $\;\lambda^{<\mu}=\lambda\;$ if $\,\mu=\aleph_0\,$
and $\;\lambda^{<\mu}\,=\,\sup\{\,\lambda^\kappa\;|\;\kappa<\mu\,\}\;$
if $\,\mu>\aleph_0\,$. Of course, 
$\;\max\{\lambda,\mu\}\leq\lambda^{<\mu}\leq\lambda^\mu\,$.
It is evident (and important)
that $\;\lambda^{<\mu}\leq\kappa<\lambda^\mu\;$
if and only if $\;\mu=\nu(\lambda,\kappa)\,$.
Therefore, referring to the spaces $\,X\in{\cal G}_\lambda\,$ in Theorem 3,
the equation $\,(w(X),|X|)=(\kappa,\kappa^{\nu(\lambda,\kappa)})\,$ 
is equivalent with $\;\lambda^{<\mu}\leq \kappa=w(X)<|X|=\lambda^\mu\,$.
\mp
For any space $\,X\,$ define $\;\hat w(X)\;$
as the smallest possible weight of a
nonempty open subspace of $\,X\,$.
Trivially, $\;\hat w(X)\leq w(X)\,$.
Clearly, any linearly ordered continuum $\,X\,$ with 
$\,\hat w(X)>\aleph_0\,$
must be totally path\-wise disconnected. 
Let us call an infinite space $\,X\,$ {\it balanced} 
if and only if $\,|U|=|X|\,$ and $\,w(U)=w(X)\,$
for every nonempty open subspace $\,U\,$ of $\,X\,$.
So $\,X\,$ is balanced if and only if 
all nonempty open subsets are equipollent 
and  $\,\hat w(X)=w(X)\,$.
(A motivation of considering these two {\it homogeneity}
properties in combination is that obviously 
any product space $\,Y^\kappa\,$ with $\,|Y|\geq 2\,$
is balanced.)
Thus the following theorem 
sharpens the statement on the family $\,{\cal F}_\kappa\,$ 
in Theorem 3 for many cardinals $\,\kappa\,$ including 
$\,\aleph_1\leq\kappa\leq \bc\,$ and all $\,\kappa\,$ satisfying 
$\;\kappa^{\aleph_0}=\kappa\,$. 
(Note that 
$\;\kappa^{\aleph_0}=\kappa\;$ if 
$\;\kappa\,\in\,\{\,2^\mu,(2^\mu)^+,(2^\mu)^{++},...\,\}\;$
for arbitrary $\,\mu\,$ in view of the Hausdorff formula [3] (5.22).)
The following theorem 
also sharpens the statement on the family $\,{\cal G}_\lambda\,$ 
if $\,\kappa=\lambda\,$ is singular.
(Note that $\,\nu(\kappa,\kappa)<\kappa\,$ 
for singular $\,\kappa\,$ due to K\"onig's theorem [3] 5.14.)
\mp
{\bf Theorem 6.} {\it If 
$\;\max\{\mu^+,\lambda^{<\mu}\}\leq\kappa\leq \lambda^\mu\;$
then there exist $\,2^\kappa\,$ 
incomparable,  balanced linearly 
ordered continua $\,X\,$ with $\,w(X)=\kappa\,$
and $\,|X|=\lambda^\mu\,$.}
\mp
Since {\it balanced} is a rather strong property, it is not surprising  
that in Theorem 6 there is less freedom of 
choosing $\,w(X)\,$ and $\,|X|\,$ than in Theorem 3.
For example, if $\,\nu(\kappa,\kappa)=\kappa\,$ 
then $\,\kappa=w(X)<|X|\,$ is excluded in Theorem 6.
(Clearly, $\;\max\{\mu^+,\lambda^{<\mu}\}\leq\kappa<\lambda^\mu\;$
is false for all $\,\lambda,\mu\,$ if $\,\nu(\kappa,\kappa)=\kappa\,$.) 
Some restrictions of the choices of $\,w(X)\,$ and $\,|X|\,$
are inevitable. Indeed, the following theorem demonstrates that 
there are certain (arbitrarily large) cardinals $\,\theta\,$ 
such that neither balanced linearly ordered compacta of size $\,\theta\,$
nor linearly ordered continua $\,X\,$ 
satisfying  $\,\hat w(X)=w(X)=|X|=\theta\,$ exist.
(Note that for every strong limit $\,\kappa\,$ the inequality 
$\,\kappa^{\aleph_0}>\kappa\,$ 
shortly says that $\,\kappa\,$ is of countable cofinality.)
\mp
{\bf Theorem 7.} {\it If $\,\kappa\,$ is a strong limit 
and $\,\kappa^{\aleph_0}>\kappa\,$ 
then a compact linearly ordered space of size 
$\,\kappa\,$ where all nonempty open sets are equipollent does not exist.
If $\,\kappa\,$ is singular
then a linearly ordered
continuum $\,X\,$ with $\,\hat w(X)=\kappa\,$ and 
$\,|X|<\kappa^{\aleph_0}\,$ does not exist.}
\mp
If $\,\kappa>\aleph_0\,$ is a strong limit 
then by Theorem 3 there exist $\,2^\kappa\,$ incomparable linearly 
ordered continua $\,X\,$ such that $\,w(X)=\kappa\,$ and $\,|X|=\theta\,$
for $\,\theta=\kappa\,$ and for $\,\theta=2^\kappa\,$, respectively.
In view of Theorem 6 and Theorem 7, 
the homogeneity property $\,\hat w(X)=w(X)\,$ can be included 
for $\,\theta=\kappa\,$ if and only if $\,\kappa^{\aleph_0}=\kappa\,$.
We claim that also for $\,\theta=2^\kappa\,$ 
the property $\,\hat w(X)=w(X)\,$ can be included 
if $\,\kappa^{\aleph_0}=\kappa\,$. 
This claim can be settled by Theorem 6 if and only if 
the uncountable strong limit $\,\kappa\,$
is singular. (If $\,\kappa\,$ is a strong limit then 
$\,\nu(\kappa,\kappa)<\kappa\,$ if and only if $\,\kappa\,$ is singular.)
In general the claim is settled by 
the following theorem which is a noteworthy 
supplement to Theorem 3 and Theorem 6.
\mp
{\bf Theorem 8.} {\it If $\,\lambda\leq\kappa\,$ and 
$\,\kappa^{\aleph_0}=\kappa\,$ then there exist $\,2^\kappa\,$ 
incomparable,  linearly 
ordered continua $\,X\,$ such that $\;\hat w(X)=w(X)=\kappa\;$
and $\,|X|=\lambda^{\nu(\lambda,\kappa)}\,$.}
\bp
{\bff 2. Proof of Theorem 1.}
\mp
We distinuguish the cases $\,\kappa=\aleph_0\,$ 
and $\,\kappa>\aleph_0\,$. Assume firstly that 
$\,\kappa>\aleph_0\,$ and 
$\,\lambda\leq\kappa\,$
and $\,\bc\leq\theta\in\{\kappa,\lambda^{\nu(\lambda,\kappa)}\}\,$.
In the following we apply Theorem 3 and Theorem 4 
which are proved in Section 5 and Section 7.
By Theorem 3 there exists a family $\,{\cal F}\,$
of incomparable totally pathwise disconnected continua 
of weight $\,\kappa\,$ and size $\,\theta\,$ such that 
$\,|{\cal F}|=2^\kappa\,$. By Theorem 4 we 
can additionally assume that all spaces in $\,{\cal F}\,$
are {\it first countable} if $\,\theta=\bc\,$.
\mp
For $\,X\in {\cal F}\,$ consider the compact product space
$\,X\times[0,1]\,$. Let $\,\Psi(X)\,$ be the quotient 
of $\,X\times[0,1]\,$ by its compact subspace $\,X\times\{1\}\,$. 
So the subspace $\,X\times \{1\}\,$ of $\,X\times[0,1]\,$
is shrinked to a point $\,p\in \Psi(X)\,$ and 
the subspace $\;X\times[0,1[\;$ of $\,X\times[0,1]\,$
is identical with the subspace 
$\,\Psi(X)\setminus\{p\}\,$ of $\,\Psi(X)\,$.
One may picture $\,\Psi(X)\,$ as a {\it cone} with the {\it basis} 
$\,X\times\{0\}\,$ 
and the {\it apex} $\,p\,$. The {\it rulings} of the cone $\,\Psi(X)\,$
are the subspaces $\;(\{x\}\times[0,1[)\cup\{p\}\;$ for $\,x\in X\,$
and they all are homeomorphic with $\,[0,1]\,$.
\mp 
Of course, $\,\Psi(X)\,$ is a compact Hausdorff space.
Since $\;|X|=\theta\geq\bc\;$ and $\,|[0,1]|=\bc\,$, we have 
$\;|\Psi(X)|=\theta\,$. 
Since $\,X\,$ is compact, by virtue of [2] 3.1.15
the filter of the neighborhoods of the compact subset 
$\,X\times\{1\}\,$ in the space $\,X\times[0,1]\,$ 
has a countable basis and hence 
$\,p\,$ has a countable local basis.
Therefore, since $\,w(X)=\kappa\,$ and $\,w([0,1[)=\aleph_0\,$,
we have $\;w(\Psi(X))=\kappa\,$. 
Furthermore, if $\,X\,$ is first countable then 
$\,\Psi(X)\,$ is first countable.
Of course, $\,\Psi(X)\,$ is pathwise connected.
Since $\,X\,$ is totally pathwise disconnected,
the apex $\,p\,$ of the cone $\,\Psi(X)\,$
is topologically labeled as the unique point $\,x\,$
where $\,\Psi(X)\setminus\{x\}\,$ has {\it infinitely many}
path components. (If $\,p\not=x\in\Psi(X)\,$
then $\,\Psi(X)\setminus\{x\}\,$ has either one or two 
path components.) On the other hand, 
the set of all $\,x\in X\,$ such that $\,\Psi(X)\setminus\{x\}\,$
remains pathwise connected is $\,X\times\{0\}\,$
and hence homeomorphic with $\,X\,$.
This enables us to verify that the $\,2^\kappa\,$ spaces 
$\;\Psi(X)\;(X\in{\cal F})\;$
are incomparable. 
\mp
Let $\;X_1,X_2\in{\cal F}\;$
and let $\,f\,$ be a homeomorphism from $\,\Psi(X_1)\,$
onto a subspace $\,S\,$ of $\,\Psi(X_2)\,$.
It is evident that $\,\Psi(X_1)\,$ is not homeomorphic 
to a subspace of $\;[0,1[\,$. 
Therefore, since $\,X_1\,$ is totally 
pathwise disconnected and $\,\Psi(X_1)\,$
is a union of $\,\kappa\,$ copies $\,C_i\,$ of $\,[0,1]\,$
with $\,C_i\cap C_j=\{p\}\,$ for $\,i\not=j\,$,
$\;f(\Psi(X_1))\subset X_2\times[0,1[\;$
is impossible and hence $\,f(p)=p\,$.
Consequently, every ruling $\,(\{x_1\}\times [0,1[)\cup\{p\}\,$
of the cone $\,\Psi(X_1)\,$ is mapped into a ruling 
$\,(\{x_2\}\times [0,1[)\cup\{p\}\,$
of the cone $\,\Psi(X_2)\,$. Therefore, since the rulings are 
homeomorphic to $\,[0,1]\,$ and $\,f(p)=p\,$, 
we can define an injective mapping $\,g\,$ from 
$\,X_1\,$ into $\,X_2\,$ and a mapping $\,h\,$
from $\,X_1\,$ into $\,[0,1[\;$ such that
$\;f(\{x\}\times [0,1[)=\{g(x)\}\times [h(x),1[\;$
for every $\,x\in X_1\,$. Thus $\;W\,:=\,\{\,(g(x),h(x))\;|\;x\in X_1\,\}\;$
is the set of all points $\,y\,$ in the pathwise connected 
space $\,S\,$ such that $\,S\setminus\{y\}\,$ remains
pathwise connected. Consequently, $\,f(X_1\times\{0\})=W\,$
and hence $\,W\,$ is a compact and connected subspace
of $\;\Psi(X_2)\setminus\{p\}\,=\,X_2\times [0,1[\;$ 
homeomorphic to $\,X_1\,$. If $\,\pi\,$ denotes the (continuous) 
projection of $\,X_2\times[0,1[\;$ onto $\,X_2\times\{0\}\,$
then $\,\pi\,$ restricted to $\,W\,$ is clearly injective.
Thus $\,\pi(W)\,$ is a compact subspace of $\,X_2\times\{0\}\,$
homeomorphic to $\,W\,$ and hence to $\,X_1\,$.
Since $\,X_2\times\{0\}\,$ is homeomorphic to $\,X_2\,$,
we conclude that $\,X_1\,$ is homeomorphic to a subspace of 
$\,X_2\,$ and hence we derive $\,X_1=X_2\,$. This 
concludes the proof of the first and the second statement
of Theorem 1 in case that $\,\kappa>\aleph_0\,$.
\mp
In order to prove the third statement
of Theorem 1 in case that $\,\kappa>\aleph_0\,$
let $\,\mu\leq\kappa\,$ and
let $\,{\cal F}_\kappa\,$ 
be a family as provided by Theorem 3.
Consider for $\,X\in{\cal F}_\kappa\,$ the cone $\,\Psi(X)\,$
and the product space $\;X^\mu\,$.
Since $\,X\,$ is a totally pathwise disconnected continuum
of weight $\,\kappa\,$ and size $\,\max\{\kappa,\bc\}\,$, 
the space $\,X^\mu\,$ is a totally pathwise disconnected continuum
of weight $\,\kappa\,$ and size $\,\kappa^\mu\,$.
For every $\,X\in{\cal F}_\kappa\,$ fix a point $\,a\in X^\mu\,$
and a point $\,b\,$ in $\;X\times\{0\}\,\subset\,\Psi(X)\,$.
(Note that, formally by the {\it Axiom of Regularity} [3], the 
sets $\,X^\mu\,$ and $\,\Psi(X)\,$ are disjoint.)
Let $\,\Phi(X)\,$ denote the quotient of the topological sum 
of the spaces $\,X^\mu\,$ and $\,\Psi(X)\,$ by the subspace 
$\;\{a,b\}\,$. (We obtain $\,\Phi(X)\,$ by
sticking together $\,X^\mu\,$ and $\,\Psi(X)\,$
at the points $\,a\,$ and $\,b\,$.) Naturally, 
$\,\Phi(X)\,$ is a continuum
of weight $\,\kappa\,$ and size $\,\kappa^\mu\,$.
\sp
Let $\,X_1,X_2\in{\cal F}_\kappa\,$ be
distinct, whence the pathwise 
connected spaces $\,\Psi(X_1)\,$ and $\,\Psi(X_2)\,$ are  
incomparable. If $\,S\,$ is a subspace           
of $\,\Phi(X_1)\,$ homeomorphic with $\,[0,1]\,$
then $\,S\,$ must be a subspace of $\,\Psi(X_1)\,$
which cannot be embedded into 
the totally pathwise disconnected
subspace $\,X_2^\mu\,$ of $\,\Phi(X_2)\,$.
Therefore, the continua 
$\,\Phi(X_1)\,$ and $\,\Phi(X_2)\,$ are incomparable as well.
This concludes the proof of all statements of Theorem 1
in case that $\,\kappa>\aleph_0\,$.
\mp
If $\,\kappa=\aleph_0\,$ then $\;\max\{\kappa,\bc\}=\bc\,$  
and $\,\lambda^{\nu(\lambda,\kappa)}=\bc\,$ and $\,\kappa^\mu=\bc\,$
and hence we are done
by the following proof of Theorem 2.
\bp
{\bff 3. Proof of Theorem 2.}
\mp
For the proof of Theorem 2 we need the following proposition.
\mp
{\bf Proposition 1.} {\it There exists a family $\;{\cal H}\;$ of 
countable sets $\,H\subset\,]{0,1}[\;$ with 
$\;\inf H=0\,$ and $\,\sup H=1\,$
such that $\;|{\cal H}|=\bc\;$
and if $\,H_1,H_2\in{\cal H}\,$ are distinct and $\,F\,$ is a finite set
then there does not exist a monotonic bijection from
$\,H_1\,$ onto $\,H_2\setminus F\,$.}
\mp
{\it Proof.} Let $\,G\,$ be the family of 
all functions from $\,\N\,$ to $\,\{3,5\}\,$, whence $\,|G|=\bc\,$.
Following [11], let $\,\zeta\,$ resp.~$\,\eta\,$ 
be the order-type of the naturally ordered set $\,\Z\,$ resp.~$\,\Q\,$.
For every $\,g\in G\,$ consider the order-type
$\;\tau[g]\,:=\,7+\zeta+g(1)+\zeta+g(2)+\zeta+g(3)+\,\cdots\,$.
Clearly, if $\,A[g]\subset\R\,$ is of order-type $\,\tau[g]\,$
for every $\,g\in G\,$ then for distinct $\,g_1,g_2\in G\,$
there is no monotonic bijection from $\,A[g_1]\,$ onto 
$\,A[g_2]\,$. Now choose for every $\,g\in G\,$
a set $\;H[g]\subset]0,1[\;$ 
of order-type 
$\;(\eta+\tau[g]+\eta+\tau[g]+\eta+\tau[g]+\,\cdots)+\zeta\;$
with $\;\inf H=0\,$ and $\,\sup H=1\,$.
Due to the {\it infinite} repetition of $\,\tau[g]\,$,
if $\,E\,$ is a {\it finite} set and $\,g_1,g_2\in G\,$ are distinct
then there is no monotonic bijection from $\,H[g_1]\,$ onto 
$\,H[g_2]\setminus E\,$. So we are done by defining
$\;{\cal H}\,=\,\{\,H[g]\;|\;g\in G\,\}\,$, {\it q.e.d.}
\mp 
Now we are ready to prove Theorem 2.
We regard the points in $\,\R^3\,$ as vectors, 
whence $\;\{\,(1-t)\cdot a+t\cdot b\;|\;0\leq t\leq 1\,\}\;$ 
is the closed straight segment which connects the points $\,a,b\in\R^3\,$.
If $\,A\subset\R^3\,$
then $\,\overline{A}\,$ denotes
the closure of $\,A\,$ in the space $\,{\Bbb R}^3\,$.  
For a pathwise connected Hausdorff space $\,X\,$
let $\,\Omega(X)\,$ denote the set of all points $\,a\in X\,$ 
where no local basis at $\,a\,$ contains only pathwise connected sets.
\mp
Let $\,{\cal H}\,$ be a family as provided by Proposition 1.
For every set $\;H\,=\,\{\,a_0,a_1,a_2,...\,\}\;$ 
in the family $\,{\cal H}\,$ with $\,a_i\not=a_j\,$ for $\,i\not=j\,$
we define a mapping 
$\;\phi=\phi_H\;$ from $\,[0,\infty[\;$ into $\,{\Bbb R}^3\,$ by 
$\;\phi(2n)\,=\,(a_n,0,0)\;\;$ and 
$\;\;\phi(2n+1)\,=\,(a_n,2^{-n},3^{-n})\;\;$  for every 
$\,n\in{\Bbb N}\cup\{0\}\,$ and by
$\;\phi(k+t)\,:=\,(1-t)\cdot \phi(k)+t\cdot \phi(k+1)\;$ 
for every integer $\,k\geq 0\,$ and every $\;t\in[0,1]\,$.
\sp
It is clear that $\,\phi\,$ is continuous. 
Furthermore $\,\phi\,$ is injective because if 
$\,\pi_n\,$ is the plane through 
the three points $\;\phi(2n),\, \phi(2n\!+\!1),\,\phi(2n\!+\!2)\;$
then $\;\pi_n\,\not=\,{\Bbb R}\times{\Bbb R}\times\{0\}\;$ and
$\;\pi_m\cap \pi_n\,=\,{\Bbb R}\!\times\!\{0\}\!\times\!\{0\}\;$
whenever $\;m,n\in{\Bbb N}\cup\{0\}\;$ and $\,m\not=n\,$.
\mp
Let $\,H\in{\cal H}\,$ and put $\;A_H\,:=\,\phi_H([0,\infty[)\,$,
whence $\,A_H\,$ is a pathwise connected subspace of $\,\R^3\,$
which is clearly not compact.
Since $\,\inf H=0\not\in H\,$ and $\,\sup H=1\not\in H\,$,
we obviously have $\;\Omega(A_H)\,=\,H\!\times\!\{0\}\!\times\!\{0\}\,$.
Since $\;A_H\subset[0,1]^3\,$, the space $\,\overline{A_H}\,$ is compact.
The space $\,\overline{A_H}\,$ is also pathwise connected
since $\;\overline{A_H}=A_H\cup I\;$
where $\;I\,:=\,[0,1]\times\{0\}\times\{0\}\,$. 
For the space $\,\overline{A_H}\,$ we obviously have 
$\;\Omega(\overline{A_H})\,=\,I\setminus\{(0,0,0),(1,0,0)\}\,$.
\mp
We claim that the $\,\bc\,$ spaces $\;\overline{A_H}\;(H\in{\cal H})\;$
are incomparable. Let $\,G,H\in {\cal H}\,$ and assume 
that $\,f\,$ is a homeomorphism from $\;\overline{A_G}\;$
onto a (compact and pathwise connected) subspace of 
$\;\overline{A_H}\,$. Then $\;f(A_G)\,=\,f\circ \phi_G([0,\infty[)\;$ 
is a pathwise connected subspace of $\,\overline{A_H}\,$ such that 
$\;f(A_G)\;$ is homeomorphic with $\,A_G\,$
and hence $\;|\Omega(f(A_G))|=\aleph_0\,$.
Naturally, $\; f\circ \phi_G\;$ is an injective continuous
mapping from $\;[0,\infty[\;$ into $\,\overline{A_H}\,$.
\mp
Let $\,g\,$ be 
any injective continuous 
mapping from $\;[0,\infty[\;$ into $\,\overline{A_H}\,$.
Then it is evident that $\;\Omega(g([0,\infty[))\subset I\,$.
Of course, $\,|\Omega(g([0,\infty[))|>\aleph_0\,$ and also 
$\,|\Omega(g([0,\infty[))|<\aleph_0\,$ is possible.
(Actually, $\,|\Omega(g([0,\infty[))|\not=\aleph_0\,$
if and only if $\,|\Omega(g([0,\infty[))|=\bc\,$
or $\,\Omega(g([0,\infty[))=\emptyset\,$.)
If $\;|\Omega(g([0,\infty[))|=\aleph_0\;$ then 
$\;g([0,\infty[)\cap I\,\subset\,\Omega(A_H)\;$
and $\;g([0,\infty[)\;$ contains 
infinitely many points of $\;(H\cap J)\times\{0\}\times\{0\}\;$ 
for all intervals $\,J=[0,\epsilon]\,$
and $\,J=[1-\epsilon,1]\,$ with $\,\epsilon>0\,$. 
As a consequence, 
$\;|\Omega(g([0,\infty[))|=\aleph_0\;$ if and only if 
$\;g([0,\infty[)=\phi_H([t,\infty[)\;$ for some $\,t\geq 0\,$. 
\mp
So we conclude that $\;f(A_G)=\phi_H([t,\infty[)\subset A_H\;$
for some $\,t\geq 0\,$.
The set $\,I\,$ is the closure of $\;\Omega(\overline{A_G})\;$
in the space $\,\overline{A_G}\,$
and $\;f(\Omega(\overline{A_G}))=\Omega(f(\overline{A_G}))\subset I\,$.
Consequently, $\,f\,$ maps $\,I\,$ onto a compact and connected
subset of $\,I\,$. We have $\,f(I)=I\,$ because
if $\,f(I)\,$ were a proper subset of $\,I\,$ then 
from $\;\phi_H([t,\infty[)\subset f(\overline{A_G})\;$
we derive that $\;f(\overline{A_G})\;$ would not be compact.
It is clear that 
$\;\Omega(\phi_H([t,\infty[))\,=\,\{\,(a_m,0,0)\;|\;2m\geq t\,\}\;$
for every $\,t\geq 0\,$. So $\,f\,$ maps 
$\;\Omega(A_G)\,=\,G\!\times\!\{0\}\!\times\!\{0\}\;$
onto a cofinite subset of
 $\;\Omega(A_H)\,=\,H\!\times\!\{0\}\!\times\!\{0\}\,$.
Of course, $\,f\,$ is essentially a monotonic permutation 
of the unit interval $\,[0,1]\,$. 
Therefore, $\,f\,$ induces a monotonic bijection 
from $\,G\,$ onto a cofinite subset of $\,H\,$.
Thus $\,G=H\,$ and this concludes the proof of Theorem 2.
\mp
{\it Remark.} With much more effort it is possible to sharpen Theorem 2.
Actually, {\it there exist $\,\bc\,$ incomparable,
pathwise connected, compact subspaces of the plane $\,\R^2\,$.}
A proof of this statement is carried out in [7].

\bp
{\bff 4. Basic constructions}
\mp
A {\it densely ordered space} is a space 
whose topology is generated by a
{\it dense} linear ordering 
(i.e.~infinitely many points lie between any two points).
It is plain that $\,d(X)=w(X)\,$ for every densely ordered space $\,X\,$.
In view of [12] 39.7, 39.8 linearly ordered continua may be characterized 
in a handy way.
\mp
{\bf Lemma 1.} {\it A space $\,X\,$ with $\,|X|\geq 2\,$
is a linearly ordered continuum 
if and only if the topology of $\,X\,$ is generated by a 
dense linear ordering of 
$\,X\,$ such that $\,X\,$ has a minimum 
and every nonempty subset of $\,X\,$
has a supremum. A space $\,X\,$ is a linearly ordered continuum 
if and only if $\,X\,$ is densely ordered and compact.}
\mp
Let $\,X\,$ be a linearly ordered continuum.  
If $\,a,b\,$ are distinct points in $\,X\,$ then
let $\,I[a,b]\,$ denote the intersection 
of all connected subsets of $\,X\,$ containing $\,a\,$ and $\,b\,$. 
Clearly, the subspace $\,I[a,b]\,$ of $\,X\,$ is a linearly ordered continuum.
We call $\,I[a,b]\,$ an {\it interval} and $\,a,b\,$ the two {\it endpoints}
of $\,I[a,b]\,$ due to the obvious fact
that if $\,\preceq\,$ is a linear ordering of $\,X\,$ generating 
the topology of $\,X\,$ and $\,a\prec b\,$ then
$\;I[a,b]\;=\;I[b,a]\;=\;\{\,x\in X\;|\;a\preceq x\preceq b\,\}\,.$
It is clear that 
$\;\hat w(X)\,=\,
\min\{\,w(I[a,b])\;|\;a,b\in X\,,\;a\not=b\,\}\,$.
Of course, the endpoints of an 
interval $\,I\,$ are purely topologically characterized as 
the two {\it noncut points} in the space $\,I\,$, i.e.~$\,x\in I\,$ is an endpoint
of $\,I\,$ if and only if $\,I\setminus\{x\}\,$ is connected.
The fact that intervals and its endpoints are topologically characterized
has two important consequences. 
Firstly, any linearly ordered continuum $\,X\,$
carries precisely {\it two} linear orderings generating the topology of 
$\,X\,$. The second consequence is the following lemma 
which is essential for our purpose and
whose prove is an easy exercise.
\mp
{\bf Lemma 2.} {\it If $\,X\,$ and $\,Y\,$ are two linearly ordered continua 
and $\,f\,$ is a bijection from $\,X\,$ onto $\,Y\,$
then $\,f\,$ is a homeomorphism if and only if
$\,f\,$ is either increasing or decreasing
referring to any linear ordering of $\,X\,$ resp.~$\,Y\,$
which generates the topology of $\,X\,$ resp.~$\,Y\,$.}
\mp
As an important consequence of Lemma 2, if $\,X_1,X_2\,$ 
are linearly ordered sets which are continua 
then $\,X_1,X_2\,$ are {\it incomparable} if and only if for 
$\,\{i,j\}=\{1,2\}\,$ and any interval $\,I[a,b]\subset X_j\,$
there is no {\it monotonic} bijection from $\,X_i\,$ onto $\,I[a,b]\,$.
\mp
Often (and throughout this paper always) it is very easy to determine
the weight of a linearly ordered space in view of the 
following lemmas. 
\mp
{\bf Lemma 3.} {\it If $\,X\,$ is a linearly ordered continuum which
has a dense subset of size $\,\kappa\,$
and contains $\,\lambda\,$ mutually disjoint intervals
then $\;\lambda\leq w(X)\leq \kappa\,$.}
\mp
{\bf Lemma 4.} {\it If $\,X\,$ is a linearly ordered continuum 
and $\,S\,$ is an infinite subspace of $\,X\,$ then 
$\;w(X)\,\geq\,w(S)\,\geq\,
|\{\,\{a,b\}\;|\;a,b\in S\,,\;a\not=b\,,\;S\cap I[a,b]=\{a,b\}\,\}|\,$.}
\mp
Lemma 3 is evident. In order to prove Lemma 4 let 
$\;{\cal A}\,=\,
\{\,\{a,b\}\;|\;a,b\in S\,,\;a\not=b\,,\;S\cap I[a,b]=\{a,b\}\,\}\;$
and let $\,{\cal B}\,$ be a basis of  $\,S\,$
and let $\,\preceq\,$ be one of the two orderings 
which generate the topology of $\,X\,$. 
For every $\,Y\in {\cal A}\,$ we may choose $\,B_Y\in{\cal B}\,$
disjoint from $\;\{\,x\in S\;|\;  \max Y\preceq x\,\}\;$
with $\,\min Y\in B_Y\,$, whence $\,B_Y\not= B_{Y'}\,$
for distinct $\,Y,Y'\in {\cal A}\,$ and therefore
$\;|{\cal B}|\geq |{\cal A}|\,$, {\it q.e.d.}
\bp
Throughout the paper, if $\,\alpha,\beta\,$ are ordinals
with $\,\alpha<\beta\,$ and $\,\beta\,$ is infinite
then we write $\,[\alpha,\beta]\,$ 
for the set of all ordinals $\,\gamma\,$ with 
$\,\alpha\leq \gamma\leq \beta\,$ and we put 
$\;[\alpha,\beta[\,:=[\alpha,\beta]\setminus\{\beta\}\,$.
(Since $\,\beta\,$ is infinite,
the writing $\,[0,1]\,$ for the Euclidean unit interval
remains unambiguous.)
For every infinite ordinal number $\,\alpha\,$
define an infinite cardinal number $\,|\alpha|\,$ by 
$\,|\alpha|:=|[0,\alpha[|\,$.
(If, as usual, $\,\alpha\,$ is identical with the set $\,[0,\alpha[\,$
then $\,|\alpha|\,$ is the {\it size of the set $\,\alpha\,$}
and the defintion of $\,|\alpha|\,$ is compatible 
with the definition $\,|S|\,$ for arbitrary sets $\,S\,$.)
\mp
For a limit ordinal $\,\gamma>0\,$ 
(which may be an infinite cardinal number) 
and a linearly ordered set $\,L\,$
the linearly ordered set $\,L^{[0,\gamma[}\,$ is the lexicographically ordered
set of all $\gamma$-sequences in $\,L\,$. So $\,a\in L^{[0,\gamma[}\,$
when $\;a\,=\,(x_\alpha)_{\alpha<\gamma}\;$ is a mapping from 
$\,[0,\gamma[\,$ 
into $\,L\,$ and $\;a=(x_\alpha)_{\alpha<\gamma}\;$ is smaller 
than $\;b=(y_\alpha)_{\alpha<\gamma}\;$ if and only if $\,a\not=b\,$
and $\,x_\beta\,$ is smaller than $\,y_\beta\,$ for the least ordinal
$\,\beta\,$ where $\,x_\beta\not=y_\beta\,$.
\mp
Assume that the linear ordering of $\,L\,$ is dense and complete
with a maximum and a minimum, whence the linearly ordered space $\,L\,$
is a continuum. Note that for every $\,\kappa\geq \bc\,$ 
there exists a linearly ordered continuum $\,L\,$ of size $\,\kappa\,$. 
Indeed, if $\,\kappa=\bc\,$ then $\,[0,1]\,$ is such a continuum.
For $\,\kappa\geq \bc\,$ consider the compact {\it long line} of weight
$\,\kappa\,$. (This generalization of the classic long line is constructed 
from $\,[0,\kappa]\,$ by placing between each
ordinal $\,\alpha<\kappa\,$
and its successor $\,\alpha+1\,$ a copy of the open unit interval 
$\;]0,1[\,$.)
\mp
Now let $\,\gamma\,$ be an infinite  cardinal number.
We claim that the linearly ordered space 
$\,L^{[0,\gamma[}\,$ is a balanced linearly ordered continuum.
Clearly the constant sequence $\,(\max L)\,$ resp.~$\,(\min L)\,$
is the maximum resp.~the minimum of $\,L^{[0,\gamma[}\,$.
Since $\,\gamma\,$ is a cardinal number,
for every limit ordinal $\,\beta<\gamma\,$
the well-ordered set $\;[\beta,\gamma[\;$ is order-isomorphic
with the well-ordered set $\;[0,\gamma[\,$.
Therefore, between any two points in $\,L^{[0,\gamma[}\,$
there lies a set of points which is
order-isomorphic to the whole set $\,L^{[0,\gamma[}\,$.
Consequently, $\,L^{[0,\gamma[}\,$ is balanced and densely ordered.
So in order to verify that $\,L^{[0,\gamma[}\,$ is a continuum
it remains to verify that every nonempty set $\,Y\subset L^{[0,\gamma[}\,$
has a supremum. 
Let $\,\sup A\,$ denote the supremum of $\,A\subset L\,$
in the linearly ordered continuum $\,L\,$ with the convention 
$\;\sup \emptyset=\min L\,$. For every ordinal $\,\alpha<\gamma\,$
let $\;\pi_\alpha\;$ denote the projection
$\;(x_\beta)_{\beta<\gamma}\mapsto x_\alpha\;$
from $\,L^{[0,\gamma[}\,$ onto $\,L\,$.
Now for $\;\emptyset\not= Y\subset L^{[0,\gamma[}\;$
define $\;y\,=\,(y_\beta)_{\beta<\gamma}\;$ recursively via
$\;y_0\,:=\,\sup\pi_0(Y)\;$ 
and 
\sp
\cl{$\;y_{\beta}\,:=\,\sup\,\pi_{\beta}
\big(Y\cap\bigcap\limits_{\alpha<\beta}\pi_\alpha^{-1}(\{y_\alpha\})\big)\;$}
\sp
for every $\,\beta\in [1,\gamma[\,$. 
A moment's reflection suffices to see that $\,y\,$
is the supremum of $\,Y\,$.
\mp
For every $\;\beta\in[1,\gamma[\;$
let $\,S_\beta\,$ denote the set of all sequences 
$\,(x_\alpha)_{\alpha<\gamma}\,$ in $\,L^{[0,\gamma[}\,$
where either $\;x_\alpha=\min L\;$ whenever
$\;\beta\leq\alpha<\gamma\;$
or $\;x_\alpha=\max L\;$ whenever
$\;\beta\leq\alpha<\gamma\,$.
Thus $\;E\,:=\,\bigcup_{\beta\in[1,\gamma[}S_\beta\;$
is the set of all 
{\it eventually constant} $\gamma$-sequences 
in $\,L\,$ where the constant value is either $\,\min L\,$ or $\,\max L\,$.
Clearly $\,E\,$ is a dense subset 
of the continuum $\,L^{[0,\gamma[}\,$.
Since $\;|S_\beta|=|L|^{|\beta|}\;$ 
for every $\,\beta\in[1,\gamma[\,$, 
we have $\;|E|=|L|^{<\gamma}\,$.
(Note that if $\,\gamma=\aleph_0\,$ then $\,|L|^{<\gamma}=|L|\,$
and that $\;|L|^{|\beta|}=|L|^\beta=|L|\;$ for $\,\beta\in\N\,$.) 
Consequently, $\;w(L^{[0,\gamma[})\leq |L|^{<\gamma}\;$ in view of Lemma 3.
\sp
We claim that $\;w(L^{[0,\gamma[})\geq|L|^{<\gamma}\,$.
This is true for $\,\gamma=\aleph_0\,$ because then $\,|L|^{<\gamma}=|L|\,$
and it is obvious that 
between any two points in $\,L^{[0,\gamma[}\,$
there lie $\,|L|\,$ mutually exclusive intervals.
(This is even true for an arbitrary infinite cardinal $\,\gamma\,$.)
So assume that $\,\gamma\,$ is an uncountable cardinal number.
By virtue of Lemma 4, for every infinite cardinal $\,\theta<\gamma\,$ 
the set $\,S_\theta\,$ is a subspace 
of $\,L^{[0,\gamma[}\,$ such that 
$\;w(L^{[0,\gamma[})\geq |S_\theta|=|L|^\theta\,$.
Consequently, $\;w(L^{[0,\gamma[})\geq |L|^{<\gamma}\,$.
This proves the claim and hence 
$\;\hat w(L^{[0,\gamma[})=w(L^{[0,\gamma[})=|L|^{<\gamma}\;$
and $\;|L^{[0,\gamma[}|=|L|^{\gamma}\;$ for every infinite cardinal 
$\,\gamma\,$.
\mp
Now we consider the special case $\,\gamma=\aleph_0\,$. 
For the sake of better reading we write $\,[0,\omega[\;$
instead of $\,[0,\aleph_0[\,$.
(The least infinite ordinal $\,\omega\,$ is identical 
with the cardinal number $\,\aleph_0\,$
due to the convention that cardinals are initial ordinals.)
The set $\,L^{[0,\omega[}\,$ is the set of all ordinary sequences 
in $\,L\,$. (The indices of a sequence are the nonnegative integers.)  
So we have $\,|L^{[0,\omega[}|=|L|^{\aleph_0}\,$ 
and $\;w(L^{[0,\omega[})=|L|^{<\omega}=|L|\,$. 
Let $\,D\,$ be a dense subset of the continuum $\,L\,$
with $\,\max L\in D\subset L\,$.
(Consequently, since $\,L\,$ is a continuum, between any two points in 
linearly ordered set $\,L\,$ there lie infinitely 
members of $\,D\,$.) 
Define a set $\,L[D]\subset L^{[0,\omega[}\,$
such that $\,(x_0,x_1,x_2,...\,)\in L^{[0,\omega[}\,$ lies in 
$\,L[D]\,$
if and only if for every index $\,n\,$
the implication 
\mp
(4.1)\qquad\qquad$\;x_n\not\in D\;\Longrightarrow\;
\forall\,m>n\,:\;x_m=\min L\;$ 
\mp
is true. In particular, $\,L[L]=L^{[0,\omega[}\,$.
Equipped with the lexicographic ordering, $\,L[D]\,$
is a compact linearly ordered space. Indeed, 
the constant sequence $\,(\min L)\,$ is the minimum of $\,L[D]\,$.
(Since $\,\max L\in D\,$,  
the constant sequence $\,(\max L)\,$ is the maximum of $\,L[D]\,$.)
In view of (4.1) the supremum 
of every nonempty set $\,Y\subset L[D]\,$
in the linearly ordered set $\,L[D]\,$ 
is equal to the supremum of $\,Y\,$ in $\,L^{[0,\omega[}\,$.
(By the way, for the {\it infimum} this is not true 
if $\;L\setminus D\;$ is infinite. Note also that 
if $\;L\setminus D\;$ is infinite then 
the continuum $\,L[D]\,$ is not a subspace of $\,L^{[0,\omega[}\,$
since the point set $\,L[D]\,$ in $\,L^{[0,\omega[}\,$
is neither connected nor compact.)
 The weight of 
$\,L[D]\,$ is $\,|D|\,$
because firstly between two distinct points in $\,L[D]\,$
there always lie $\,|D|\,$ copies of $\,L[D]\,$
and secondly
the eventually constant 
sequences $\;(d_0,d_1,...,d_n,a,a,a,...\,)\;$ with $\,d_i\in D\,$
and $\,a=\min L\,$
form a dense subset of $\,L[D]\,$ of size $\,|D|\,$.
In particular, $\,L[D]\,$ is a balanced linearly ordered continuum
and $\;|L[D]|\,=\,\max\{|L|,|D|^{\aleph_0}\}\;$ 
and $\,w(L[D])=|D|\,$.
It is plain that if $\,L\,$ is first countable 
then $\,L[D]\,$ is first countable.
More generally, it is plain that $\,\chi(L[D])\leq \chi(L)\,$.
\mp
We conclude this section with a proposition which is essential
for the proof of Theorem 6.
\mp
{\bf Proposition 2.} {\it If $\;\lambda^{<\mu}\leq\kappa\leq\lambda^\mu\;$
then there exists 
a balanced linearly ordered continuum $\,K\,$ such that 
$\,w(K)=\kappa\,$ and $\,|K|=\lambda^\mu\,$ and
$\,\chi(K)\leq \max\{\tilde\lambda,\mu\}\,$
with $\,\tilde\lambda=\aleph_0\,$ for $\,\lambda\leq \bc\,$
and $\,\tilde\lambda=\lambda\,$ for $\,\lambda>\bc\,$.}
\mp
{\it Proof.} Assume $\;\lambda^{<\mu}\leq\kappa\leq \lambda^\mu\,$.
We distinguish the cases 
$\,\lambda^{<\mu}<\bc\,$ and $\,\lambda^{<\mu}\geq \bc\,$.
Assume firstly that $\,\lambda^{<\mu}<\bc\,$.
Then $\,\mu=\aleph_0\,$ and $\,\lambda<\bc\,$
and hence $\,\lambda^\mu=\bc\,$.  
Thus we are done by taking $\,[0,1]\,$
for the desired continuum $\,K\,$
if $\,\kappa=\aleph_0\,$ and by putting 
$\;K=L[D]\;$ where $\,L=[0,1]\,$ and 
and $\,D\subset [0,1]\,$ is dense and $\,1\in D\,$
and $\,|D|=\kappa\,$ if $\,\kappa>\aleph_0\,$. 
(In both cases $\,K\,$ is first countable.)
\sp
Now assume that $\,\lambda^{<\mu}\geq \bc\,$.
We distinguish the cases $\,\lambda\leq \bc\,$ and $\,\lambda>\bc\,$.
Assume firstly that $\,\lambda\leq \bc\,$.  
Consider the balanced linearly ordered continuum 
$\,L=[0,1]^{[0,\mu[}\,$ of weight $\,\bc^{<\mu}\,$
and size $\,\bc^\mu\,$. 
Since $\,\theta^\mu=c^\mu\,$ whenever $\,\aleph_0\leq\theta\leq \bc\,$
and since $\,\lambda\leq c\leq \lambda^{<\mu}\,$,
we have $\,\bc^\mu=\lambda^\mu\,$
and $\,\bc^{<\mu}=\lambda^{<\mu}\,$.
(So we are already done by putting $\,K=L=[0,1]^{[0,\mu[}\,$
if $\,\kappa=\lambda^{<\mu}\,$.) 
Since $\,d(L)=w(L)\leq\kappa\leq |L|\,$,
we can choose a dense set $\,D\subset L\,$
with $\,\max L\in D\,$ and $\,|D|=\kappa\,$. 
Consider $\;L[D]\subset L^{[0,\omega[}\;$ and
put $\;K\,:=\, L[D]\,$. 
Then $\,K\,$ is a balanced linearly ordered continuum 
with $\,w(K)=\kappa\,$ and 
$\,|K|=\max\{\lambda^\mu,\kappa^{\aleph_0}\}=\lambda^\mu\,$.
We have $\;\chi(K)\leq \chi([0,1]^{[0,\mu[})\leq \mu\;$
since $\,[0,1]\,$ is first countable and $\,|[0,\mu[|=\mu\,$.
\sp
Assume secondly that $\,\lambda>\bc\,$, whence  $\,\lambda^{<\mu}>\bc\,$.
Fix any linearly ordered continuum $\,\tilde L\,$ 
of size $\,\lambda\,$. Then $\,L=\tilde L^{[0,\mu[}\,$
is a balanced linearly ordered continuum 
of weight $\,\lambda^{<\mu}\,$ and size $\,\lambda^\mu\,$.
As above we choose a dense set $\,D\subset L\,$
with $\,\max L\in D\,$ and $\,|D|=\kappa\,$
and put $\;K\,:=\, L[D]\,$. Then $\,K\,$
is a balanced linearly ordered continuum 
of weight $\,\kappa\,$ and size $\,\lambda^\mu\,$.
And we have $\;\chi(K)\leq \chi(\tilde L^{[0,\mu[})\,$. 
It is plain that 
$\;\chi(\tilde L^{[0,\mu[})\leq\max \{\chi(\tilde L),\mu\}\,$.
Since $\,|\tilde L|=\lambda\,$ and 
$\,\chi(X)\leq |X|\,$ for every linearly ordered space $\,X\,$,
we derive $\;\chi(K)\leq\max\{\lambda,\mu\}\,$, {\it q.e.d.}

\bp
{\bff 5. Proofs of Theorem 3 and Theorem 6}
\mp
Let $\,\Sigma_\kappa\,$ denote the family of all 
functions from $\;[0,\kappa[\;$ into the set $\,\{0,1\}\,$.
So the members of $\,\Sigma_\kappa\,$ are $\kappa$-sequences
whose values are either $\,0\,$ or $\,1\,$.
Since $\;|[0,\kappa[|=\kappa\,$, we have
$\;|\Sigma_\kappa|=2^\kappa\,$.
Define an equivalence relation $\,\sim\,$ on $\,\Sigma_\kappa\,$
in the following way.
For $\,\sigma_1,\sigma_2\in\Sigma_\kappa\,$
define $\,\sigma_1\sim\sigma_2\,$
if and only if there is an ordinal 
$\,\gamma<\kappa\,$ such that $\;\sigma_1(\alpha)=\sigma_2(\alpha)\;$
for every $\;\alpha\in[\gamma,\kappa[\;$.
Let us call $\,\sigma_1,\sigma_2\in\Sigma_\kappa\,$
{\it similar} when $\,\sigma_1\sim\sigma_2\,$.
\mp
(5.1)$\;$ {\it The family $\,\Sigma_\kappa\,$ contains $\,2^\kappa\,$
pairwise non-similar $\kappa$-sequences.}       
\mp
Let $\,{\cal E}_\kappa\,$ denote the partition of $\,\Sigma_\kappa\,$
into equivalence classes with respect to $\,\sim\,$.
(5.1) is verified by proving that $\,|{\cal E}_\kappa|=2^\kappa\,$.
First of all we note that if $\,\mu<\kappa\,$
then \hbox{$\,|{\cal E}_\kappa|\geq 2^\mu\,$.}
Indeed, for every $\mu$-sequence $\,g\in\Sigma_\mu\,$
put $\;\sigma[g](\mu\cdot\beta+\alpha)=g(\alpha)\;$ 
whenever $\;\alpha\in[0,\mu[\;$ and $\,\beta\in[0,\kappa[\,$.
Referring to ordinal arithmetics, $\,\sigma[g]\,$
is a well-defined $\kappa$-sequence in $\,\Sigma_\kappa\,$
and we always have $\;\sigma[g_1]\not\sim\sigma[g_2]\;$
for distinct $\,g_1,g_2\in\Sigma_\mu\,$.
Consequently, $\;|{\cal E}_\kappa|\geq|\Sigma_\mu|=2^\mu\,$. 
\mp
Since $\,|{\cal E}_\kappa|\geq 2^\mu\,$ for all $\,\mu<\kappa\,$, 
we have $\,|{\cal E}_\kappa|\geq 2^{<\kappa}\,$.
Therefore $\,|{\cal E}_\kappa|=2^\kappa\,$ 
in case that $\,2^{<\kappa}=2^\kappa\,$.
In order to finish the proof of (5.1) it remains to verify 
$\,|{\cal E}_\kappa|=2^\kappa\,$ in case that $\,2^{<\kappa}<2^\kappa\,$.
This is obviously accomplished by verifying that $\,|E|\leq 2^{<\kappa}\,$
for every equivalence class $\,E\in{\cal E}_\kappa\,$.
\mp
If $\,\kappa=\aleph_0\,$ (and hence $\,2^{<\kappa}=\aleph_0\,$)
then it is clear that $\,|E|=\aleph_0\,$ 
for every $\,E\in{\cal E}_\kappa\,$.
If $\,\kappa>\aleph_0\,$ then for every 
$\,E\in{\cal E}_\kappa\,$ we have
$\;|E|\,\leq\,\sum_{\omega\leq\gamma<\kappa}|\Sigma_{|\gamma|}|
=\,\sum_{\omega\leq\gamma<\kappa}2^{|\gamma|}\,=\,2^{<\kappa}\;$
and this concludes the proof of (5.1).
\mp\sp
If $\,\gamma<\kappa\,$ is an ordinal number then, naturally,
the canonically well-ordered sets $\,[0,\kappa[\;$
and $\,[\gamma,\kappa[\;$ are order-isomorphic
and there is precisely one order isomorphism 
from $\,[0,\kappa[\;$ onto $\,[\gamma,\kappa[\,$.
Moreover, the following is true.
\mp
(5.2)$\;$ {\it If $\,\kappa>\aleph_0\,$ and $\,\gamma\in[0,\kappa[\;$
and $\,\varphi\,$ is an increasing bijection from 
$\,[0,\kappa[\;$ onto $\,[\gamma,\kappa[\;$ then there
is an ordinal $\,\beta\in[\gamma,\kappa[\;$ such that 
$\;\varphi(\alpha)=\alpha\;$ for every $\;\alpha\in[\beta,\kappa[\,$.}
\mp
{\it Proof.} If $\,\gamma=0\,$ then $\,\varphi\,$ 
is the identity and we are done.
Assume $\,\gamma>0\,$ and
consider the (countable) set 
$\;A\,=\,\{\,\phi_n\;|\;n<\omega\,\}\;$
where $\,\phi_0=0\,$ and $\,\phi_1=\varphi(0)\,$
and $\,\phi_2=\varphi(\varphi(0))\,$ and 
$\,\phi_3=\varphi(\varphi(\varphi(0)))\,$ and so on, 
whence $\,\phi_0<\phi_1<\phi_2<\phi_3<\cdots\,$.
Then $\,\beta=\sup A\,$ is a limit ordinal in $\,[0,\kappa]\,$
and since $\,\sup \varphi(A)=\varphi(\sup A)\,$ we have 
$\,\varphi(\beta)=\beta\,$. We claim that $\,\beta<\kappa\,$.
(This is clearly true if $\,\kappa\,$ is of uncountable cofinality.)
Put $\;\theta=|[0,\varphi(0)[|\,$, whence $\,\theta<\kappa\,$.
Then $\;|[\phi_n,\phi_{n+1}[|=\theta\;$ 
for every $\,n\in\N\,$. Clearly, 
$\;[0,\beta[\;=\,\bigcup_{n=0}^\infty [\phi_n,\phi_{n+1}[\,$.
Consequently, $\;|[0,\beta[|=\max\{\aleph_0,\theta\}<\kappa\;$ 
and hence $\,\beta<\kappa\,$.
Finally, since $\,\varphi(\beta)=\beta\,$,
it is clear that $\,\varphi\,$ restricted to $\,[\beta,\kappa[\;$
must be the identity, {\it q.e.d.}
\mp\sp
If $\,X\,$ is a Hausdorff space and $\,x\in X\,$ then, as usual, 
$\,\chi(x,X)\,$ denotes the least possible size of a 
neighborhood basis at the point $\,x\,$ in the space $\,X\,$.
Consequently, $\;\chi(X)\,=\,\sup\,\{\,\chi(x,X)\;|\;x\in X\,\}\,$.
(Note that this supremum is not necessarily a maximum. For example
if $\,X\,$ is the compact linearly ordered space $\,[0,\kappa]\,$
for some singular $\,\kappa\,$ then $\,\chi(X)=\kappa\,$ 
but $\,\chi(x,X)<\kappa\,$ for every $\,x\in X\,$.)
In order to prove Theorem 3 and Theorem 6 we will apply 
the following theorem.
\mp
{\bf Theorem 10.} {\it Let $\,K\,$ be a linearly ordered continuum 
such that $\,\chi(K)<w(K)\leq\kappa\,$.
Then there exists a family $\,{\cal L}_\kappa(K)\,$
of incomparable linearly ordered continua 
such that $\;|{\cal L}_\kappa(K)|=2^\kappa\;$
and $\,w(X)=\kappa\,$ and $\,\hat w(X)=\hat w(K)\,$ 
and $\,|X|=\max\{\kappa,|K|\}\,$ for every $\,X\in{\cal L}_\kappa(K)\,$
and if $\,K\,$ is balanced and $\,w(K)=\kappa\,$
then every space in $\,{\cal L}_\kappa(K)\,$ is balanced.}
\mp
{\it Proof.} Let $\,K\,$ be a linearly ordered set which is a continuum 
in the order topology 
such that $\,\theta:=\chi(K)<w(K)\leq\kappa\,$.
Then $\,\theta^+\,$ is a regular cardinal 
and $\,\theta^+\leq w(K)\leq\kappa\,$.
Let a modification $\,M\,$ of the compact long line be constructed from 
$\,[0,\theta^+]\,$ by placing between each
ordinal $\;\alpha<\theta^+\;$
and its successor $\,\alpha+1\,$ a copy of the densely ordered set
$\;K\setminus\{e_1,e_2\}\;$
where $\,e_1, e_2\,$ are the two endpoints of $\,K\,$.
Then $\,M\,$ is a compact, densely ordered set 
with minimum $\,0\,$ and maximum $\,\theta^+\,$
and $\,w(M)=\kappa\,$ and $\,|M|=|K|\,$. 
Since $\,\theta^+\,$ is a regular cardinal, 
the maximum $\,\theta^+\,$ of $\,M\,$ is the unique point $\,x\,$ 
in $\,M\,$ where $\;\chi(x,M)>\theta\,$.
Let $\;Y\,=\,M\setminus\{0,\theta^+\}\;$ and
let $\,Y^*\,$ be the set $\,Y\,$
equipped with the {\it backwards linear ordering} of $\,Y\,$.
(If $\,x,y\in Y\,$ then $\,x\,$ is smaller than $\,y\,$
with respect to the original ordering of $\,Y\,$ 
if and only if $\,y\,$ is smaller than $\,x\,$
with respect to the backwards ordering.)
\mp
Let $\,\Gamma\,$ denote the set of all limit ordinals
in the set $\;[\omega,\kappa[\,$
and let $\,\Sigma\,$ be the family of all functions from 
$\,\Gamma\,$ into the set $\,\{0,1\}\,$.
Naturally, $\;|\Gamma|=\kappa\;$ and hence
$\;|\Sigma|=2^\kappa\,$.
Naturally, $\,\Gamma\,$ is well-ordered 
and (hence) order-isomorphic with the 
canonically ordered set $\,[0,\kappa[\,$.
Therefore we can speak of {\it similar} 
and {\it non-similar} members of $\,\Sigma\,$
referring to the definition of the equivalencec relation $\,\sim\,$
defined above.
\mp
For each $\;\sigma\in\Sigma\;$
let $\,X_\sigma\,$ be constructed from 
$\,[0,\kappa]\,$ by placing between each ordinal 
$\,\alpha<\kappa\,$
and its ordinal successor $\,\alpha+1\,$ a copy of either the ordered
set $\,Y\,$ or the ordered set $\,Y^*\,$ in the following way.
If $\;\alpha\not\in \Gamma\;$ 
then we place between $\,\alpha\,$ and $\,\alpha+1\,$ a copy of $\,Y\,$.
If $\,\gamma\in\Gamma\,$ then between $\,\gamma\,$ and $\,\gamma+1\,$
we place a copy of $\,Y\,$
when $\;\sigma(\gamma)=0\;$ and 
a copy of $\,Y^*\,$ when $\;\sigma(\gamma)=1\,$.
In this way the set $\,X_\sigma\,$ becomes linearly ordered
with minimum $\,0\,$ and maximum $\,\kappa\,$ and
the linear ordering of $\,X_\sigma\,$ 
restricted to the subset 
$\,[0,\kappa]\,$ of $\,X_\sigma\,$ 
is the canonical well-ordering of $\,[0,\kappa]\,$.
Clearly, with respect to the order topology,
$\,X_\sigma\,$ is a linearly ordered continuum and
$\;|X_\sigma|=\max\{\kappa,|M|\}=\max\{\kappa,|K|\}\;$
and  $\,\hat w(X_\sigma)=\hat w(K)\,$.
By virtue of Lemma 3 we have $\,w(X_\sigma)=\kappa\,$.
It is evident that if $\,w(K)=\kappa\,$ and 
(and hence $\,|X_\sigma|=|K|\,$) and $\,K\,$
is balanced then the space $\,X_\sigma\,$ is balanced.
\mp
The two endpoints $\,0,\kappa\,$ of the linearly ordered
continuum $\,X_\sigma\,$
can be distinguished purely topologically.
For the point $\,0\,$
has a compact neighborhood which is a space of 
character not greater than $\,\theta\,$, while 
the point $\,\kappa\,$ has no such neighborhood.
(If the cofinality of $\,\kappa\,$ is greater than $\,\theta\,$
then there is also the distinction that 
$\,\chi(0,X_\sigma)\leq\theta\,$ and
$\,\chi(\kappa,X_\sigma)>\theta\,$.)
Therefore, the linear ordering $\,\preceq_\sigma\,$ which generates 
the topology of the continuum $\,X_\sigma\,$ 
can be chosen so that $\,0\,$ is the minimum 
of $\,X_\sigma\,$ and $\,\kappa\,$ is the maximum.
Then $\,\preceq_\sigma\,$ is completely determined 
by the topology of the continuum $\,X_\sigma\,$.
The subset $\,\Gamma\,$ of $\,X_\sigma\,$
is also completely determined 
by the topology of $\,X_\sigma\,$.
Indeed, let $\,W_\sigma\,$ denote 
the set of all points $\,x\,$ in the space 
$\;X_\sigma\setminus\{\kappa\}\;$
such that $\;\chi(x,X_\sigma\setminus\{\kappa\})>\theta\,$.
Then we obviously have   
\sp
\cl{$W_\sigma\;=\;[0,\kappa]\setminus(\{0,\kappa\}\cup
\{\,\gamma+1\;\,|\,\;\gamma\in\Gamma\;\land\;\sigma(\gamma)=1\,\})$}
\sp
and hence $\,(W_\sigma,\preceq_\sigma)\,$ 
is order-isomorphic with $\,[0,\kappa[\;$
and $\,\Gamma\,$ is 
the set of all limit points in the subspace $\,W_\sigma\,$
of the continuum $\,X_\sigma\,$. 
\sp
Let $\,\sigma\in\Sigma\,$.
For every $\,\gamma\in\Gamma\,$ 
let $\,\xi(\gamma)\,$ be an arbitrary point in $\,X_\sigma\,$
with $\,\xi\succ\gamma\,$. (For example, let 
$\,\xi(\gamma)\,$ be the ordinal number $\,\gamma+2\,$ which lies in 
$\;W_\sigma\setminus\Gamma\,$.)
For $\,\gamma\in\Gamma\,$ consider the 
interval $\,I[\gamma,\xi(\gamma)]\,$ in the space $\,X_\sigma\,$.
Then we obviously have 
$\;\chi(\gamma,I[\gamma,\xi(\gamma)])>\theta\;$
if and only if $\;\sigma(\gamma)=1\,$.
Therefore, every $\,\sigma\in\Sigma\,$ 
is completely determined by the topology of $\,X_\sigma\,$.
\sp
We claim that if $\,\sigma,\sigma'\in \Sigma\,$
are not similar then the spaces $\,X_\sigma,X_{\sigma'}\,$
are incomparable. Then we are done by virtue of (5.1).
Let $\,\sigma,\sigma'\in \Sigma\,$ and let 
$\,f\,$ be a homeomorphism from $\,X_\sigma\,$ 
onto a subspace $\,S\,$ of $\,X_{\sigma'}\,$.
Of course, $\,S\,$ is a continuum and hence 
$\,S\,$ is an interval in $\,X_{\sigma'}\,$
and hence $\;S\,=\,\{\,x\;|\;u\preceq_{\sigma'}x\preceq_{\sigma'}v\,\}\;$
for some $\,u,v\in X_{\sigma'}\,$ with $\,u\prec_{\sigma'} v\,$.
Since in $\,X_{\sigma}\,$ and in $\,X_{\sigma'}\,$ the point
$\,\kappa\,$ is uniquely topologically labeled,
$\,f(\kappa)=\kappa=v\,$ and hence
$\,f\,$ is a strictly increasing mapping 
from $\,X_\sigma\,$ onto $\,S\subset X_{\sigma'}\,$. 
Since $\,f\,$ is a homeomorphism from $\,X_\sigma\,$ onto $\,S\,$,
we must have $\;f(W_\sigma)\subset W_{\sigma'}\;$ 
and hence $\;f(\Gamma)=\Gamma\cap S=\Gamma\cap[\eta,\kappa[\;$
for some $\,\eta\in \Gamma\,$. 
Since $\,\Gamma\,$ is order-isomorphic with $\,[0,\kappa[\,$,
by virtue of (5.2) we can find 
an ordinal $\,\gamma\,$ in $\;\Gamma\cap[\eta,\kappa[\;$
such that  $\,f(\alpha)=\alpha\,$
for every  $\;\alpha\,\in\,\Gamma\cap[\eta,\kappa[\,$.
Therefore, since $\,\sigma\,$ and $\,\sigma'\,$ 
are completely determined by the topologies of $\,X_\sigma\,$
and $\,X_{\sigma'}\,$, we must have $\;\sigma(\alpha)=\sigma'(\alpha)\;$
whenever $\;\alpha\in \Gamma\;$ and $\,\alpha\geq\eta\,$.
This proves the claim and concludes 
the proof of Theorem 10, {\it q.e.d.} 
\mp\sp
{\it Proof of Theorem 3.} Assume $\,\kappa>\aleph_0\,$.
Referring to Section 4, 
choose any dense $\,D\subset [0,1]\,$ such that $\,1\in D\,$
and $\,|D|=\aleph_1\,$ and put $\,K=L[D]\,$ for $\,L=[0,1]\,$.
Then $\,K\,$ is a {\it first countable}
linearly ordered continuum with $\,w(K)=\hat w(K)=\aleph_1\,$, whence
$\,|K|=\bc\,$.
Consequently, by applying Theorem 10,
there exists a family $\,{\cal F}_\kappa={\cal L}_\kappa(K)\,$ 
of incomparable linearly ordered continua such that 
$\,|{\cal F}_\kappa|=2^\kappa\,$
and $\,w(X)=\kappa\,$ and $\,\hat w(X)=\aleph_1\,$
and $\,|X|=\max\{\kappa,\bc\}\,$ for every $\,X\in{\cal F}_\kappa\,$.
In particular,  every space $\,X\,$ in 
$\,{\cal F}_\kappa\,$ is totally pathwise disconnected. 
\sp
In order to conclude the proof of Theorem 3 by 
constructing also families $\,{\cal G}_\lambda\,$ as desired,
there is nothing to show if 
$\,\lambda\leq\kappa<\bc\,$ because then 
$\,\nu(\lambda,\kappa)=\aleph_0\,$ and hence
$\,\max\{\kappa,\bc\}=c=\lambda^{\nu(\lambda,\kappa)}\,$
and thus we can put $\,{\cal G}_\lambda={\cal F}_\kappa\,$. 
So assume  $\,\kappa\geq \bc\,$ and $\,\lambda\leq\kappa\,$.
Since $\;\lambda^{<\mu}\leq\kappa<\lambda^\mu\;$
if and only if $\;\mu=\nu(\lambda,\kappa)\,$,
by applying Proposition 2  we can fix
a linearly ordered continuum $\,G\,$ 
such that $\,\hat w(G)=w(G)=\kappa\,$ 
and $\,|I|=\lambda^{\nu(\lambda,\kappa)}\,$ for each interval $\,I\,$ in 
$\,G\,$. (Since $\,\hat w(G)\geq \aleph_1\,$, $\,G\,$ 
is totally pathwise disconnected.)
Define $\;{\cal G}_\lambda\,=\,\{\,G\lor X\;|\;X\in{\cal F}_\kappa\,\}\;$
where $\,G\lor X\,$ is 
obtained by sticking together 
the linearly ordered continua
$\,G\,$ and $\,X\,$ at one common endpoint.
More precisely, $\;G\lor X\,=\tilde G\cup\tilde X\;$
for some linearly ordered copies $\,\tilde G\,$ resp.~$\,\tilde X\,$ 
of $\,G\,$ resp.~$\,X\,$ with $\,|\tilde G\cap\tilde X|=1\,$ 
such that $\,\tilde G\,$ and $\,\tilde X\,$ are intervals in 
$\;G\lor X\,$. Clearly, each space in $\,{\cal G}_\lambda\,$
is a totally pathwise disconnected
linearly ordered continuum of weight $\,\kappa\,$ 
and size $\,|G|=\lambda^{\nu(\lambda,\kappa)}\,$.
For distinct and hence incomparable 
spaces $\,X_1,X_2\in{\cal F}_\kappa\,$ the spaces
$\,G\lor X_1\,$ and $\,G\lor X_2\,$ are incomparable as well.
Indeed, let $\;X_1,X_2\in{\cal F}_\kappa\;$ 
and let $\,f\,$ be a homeomorphism from 
$\,G\lor X_1\,$ onto a subspace of $\,G\lor X_2\,$.
Since every neighborhood of $\,a\in G\,$ in the space  
$\,G\lor X_1\,$ has size $\,|G|\,$ and 
every point $\,b\,$ in $\,(G\lor X_i)\setminus G\,$
has a neighborhood of size smaller than $\,|G|\,$,
we must have $\;f(G)\subset G\,$. 
Hence $\;f(X_1)\subset X_2\;$ and this 
means that $\,X_1\,$ and $\,X_2\,$ are not incomparable, 
whence $\,X_1=X_2\,$, {\it q.e.d.}
\bp
{\it Proof of Theorem 6.} Theorem 6 is an immediate consequence of Theorem 10
and Proposition 2. Indeed, 
if $\;\max\{\mu^+,\lambda^{<\mu}\}\leq\kappa\leq \lambda^\mu\;$
then by virtue of Proposition 2 we can find a
balanced linearly ordered continuum $\,K\,$
such that $\;w(K)=\kappa\;$ 
and $\;|K|=\lambda^\mu\;$ and $\;\chi(K)\leq\mu <\kappa\,$. 
\bp
There is a consequence of Theorem 6 worth mentioning. 
\mp
{\bf Corollary 1.} {\it If $\,\kappa^{\aleph_0}\leq\lambda\,$ and 
a continuum of size $\,\lambda\,$ and weight smaller than $\,\kappa\,$
exists then there exist $\,2^\kappa\,$ 
incomparable continua of weight $\,\kappa\,$ and size $\,\lambda\,$.}    
\mp
{\it Proof.} Let $\,Y\,$ be a continuum with with $\,w(Y)<\kappa\,$ and
$\,\kappa^{\aleph_0}\leq |Y|=\lambda\,$. 
We may assume that $\,Y\,$ has no cut points.
(Otherwise we replace $\,Y\,$ with $\,Y\times Y\,$.
Of course, $\,w(Y\times Y)=w(Y)\,$ and $\,|Y\times Y|=|Y|\,$.) 
Apply Theorem~6 with $\,\mu=\aleph_0\,$ and $\,\lambda=\kappa\,$
and let $\,{\cal F}\,$ be a family of incomparable, linearly ordered 
continua $\,X\,$ such that $\,|{\cal F}|=2^\kappa\,$
and $\,\hat w(X)=w(X)=\kappa\,$ and $\,|X|=\kappa^{\aleph_0}\,$
for each $\,X\in{\cal F}\,$. 
Now fix a point $\,y\in Y\,$ and let $\,a_X\,$ be one of the two noncut points 
of $\,X\,$ for each $\,X\in{\cal F}\,$. For every $\,X\in{\cal F}\,$
let $\,Y\lor X\,$ denote the continuum obtained by attaching 
$\,X\,$ to $\,Y\,$ by identifying the point $\,a_X\in X\,$
with the point $\,y\in Y\,$. 
Let $\,X_1,X_2\in{\cal F}\,$ and assume that $\,f\,$ is a
continuous injection from $\,Y\lor X_1\,$ to $\,Y\lor X_2\,$.
We have $\,f(y)=y\,$ because firstly $\,f(y)\,$ cannot lie in 
$\,X_2\setminus\{y\}\,$ since   
the point $\,y\,$ is obviously 
the unique point $\,x\,$ in the continuum $\,Y\lor X_1\,$
such that every neighborhood of $\,x\,$ contains cut points and noncut points. 
And secondly, $\,f(y)\,$ cannot lie in 
$\,Y\setminus\{y\}\,$ since 
$\,w(U)>w(Y)\,$ for every nonempty open $\,U\subset X_1\,$ containing 
$\,y\,$. Trivially, for $\,i\in\{1,2\}\,$ the space 
$\,(Y\lor X_i)\setminus\{y\}\,$ 
has the two components $\,Y\setminus\{y\}\,$ and $\,X_i\setminus\{y\}\,$.
Consequently, either $\,f(X_1\setminus\{y\})\subset X_2\setminus\{y\}\,$ 
or $\,f(X_1\setminus\{y\})\subset Y\setminus\{y\}\,$.
But $\,f(X_1\setminus\{y\})\subset Y\setminus\{y\}\,$ is impossible
due to $\,w(f(X_1\setminus\{y\}))=\kappa>w(Y)\,$. 
So $\,f(X_1)\subset X_2\,$ 
and hence $\,f\,$ (restricted to the continuum $X_1$)
is a homeomorphism from $\,X_1\,$
onto a subspace of $\,X_2\,$ and hence $\,X_1=X_2\,$.
Therefore, the family $\,\{\,Y\lor X\;|\;X\in{\cal F}\,\}\,$
consists of $\,2^\kappa\,$ incomparable continua 
of weight $\,\kappa=\max\{w(Y),\kappa\}\,$ 
and size $\,\lambda=\max\{|Y|,\kappa^{\aleph_0}\}\,$, {\it q.e.d.}
\mp\sp
{\it Remark.} In comparison with Corollary 1, 
{\it if $\,Y\,$ is a continuum and $\,w(Y)\leq\kappa\leq |Y|\,$
then there exist $\,2^\kappa\,$ 
mutually non-homeomorphic continua of weight $\,\kappa\,$
and size $\,|Y|\,$.} This enumeration theorem is evident 
in view of the proofs of [9] Propositions 2 and 3.
\bp\sp
{\bff 6. Proof of Theorem 5}
\mp
The main idea in the proof of Theorem 5 will lead us to a proof 
of Theorem 4. Therefore we prove Theorem 5 in this section 
and Theorem 4 in the next section.
First of all notice the following well-know fact (cf.~[2] 6.3.8.c).
\sp
(6.1)\quad{\it  Every separable linearly ordered continuum is 
homeomorphic to the space $\,[0,1]\,$.}
\mp
If $\,X\,$ is a linearly ordered continuum then the {\it path components}
of $\,X\,$ are not necessarily compact subspaces of $\,X\,$.
For example, the classical compact and connected {\it long line} of weight 
$\,\aleph_1\,$  has precisely 
two path components, one of them is a singleton and
the other one is not compact. 
The following lemma is essential for the proof of 
Theorem 5.
\mp
{\bf Lemma 5.} {\it If $\,X\,$ is a first countable 
linearly ordered continuum then 
every path component of $\,X\,$ is a compact subspace of $\,X\,$
and hence either a singleton or
homeomorphic to the Euclidean space $\,[0,1]\,$.}
\mp
{\it Proof.} Let $\,X\,$ be a  linearly ordered set, compact and 
connected and first countable in the order topology.
Let $\,A\,$ be a non-singleton 
path component of $\,X\,$ 
and let $\,s\,$ be the supremum of $\,A\,$. Since $\,s\,$
has a countable neighborhood base, 
$\,s\,$ is the limit of some strictly increasing sequence 
$\;a_1,a_2,a_3,...\;$ of points in $\,A\,$.
Since $\;I[a_n,a_{n+1}]\subset A\;$ 
is homeomorphic to $\,[0,1]\,$ for every $\,n\,$,
it is plain to define a continuous mapping $\,f\,$
from $\,[0,1]\,$ into $\,X\,$ with $\,f(0)=a_1\,$ and $\,f(1)=s\,$,
whence $\,s\in A\,$. Similarly, the infimum of $\,A\,$ also lies in $\,A\,$,
{\it q.e.d.}
\sp\mp
Let $\,{\cal D}\,$ be the family 
of all subsets $\,D\,$ of $\,[0,1]\,$ such that 
$\;|D\cap[x,y]|=\aleph_1\;$ whenever $\;0\leq x<y\leq 1\,$
and for $\,D\in{\cal D}\,$ consider the set 
\sp
\cl{$C[D]\;\,:=\,\;([0,1]\setminus D)\!\times\!\{0\}\;\cup\;
D\!\times\![0,1]$}
\sp
equipped with the lexicographic ordering.
One can say that $\,C[D]\,$ is constructed from the naturally ordered 
set $\,[0,1]\,$  by replacing each 
point in $\,D\,$ with a copy of $\,[0,1]\,$.
Obviously the space $\,C[D]\,$ is a 
first countable linearly ordered continuum. 
And by applying Lemma~3 we see that $\,w(C[D])=\aleph_1\,$.
The {\it non-singleton path components}
of the space $\,C[D]\,$ are precisely the compact sets 
$\;\{d\}\!\times\![0,1]\;$ with $\,d\in D\,$.
Because $\;\{d\}\!\times\![0,1]\;$ is 
trivially pathwise connected for every $\,d\in D\,$ and 
if $\,x,y\in C[D]\,$ 
and $\,x\not=y\,$ and $\;\{x,y\}\not\subset\{d\}\!\times\![0,1]\;$
for every $\,d\in D\,$
then in the linearly ordered set $\,C[D]\,$ we find uncountably many 
mutually disjoint copies of $\,[0,1]\,$ lying between $\,x\,$ and $\,y\,$ 
and hence there does not exist a continuous mapping $\,f\,$ from 
$\,[0,1]\,$ into $\,C[D]\,$ with $\,f(0)=x\,$ and $\,f(1)=y\,$.
The union of the singleton path components of $\,C[D]\,$ 
is $\;([0,1]\setminus D)\times\{0\}\;$ and hence a
separable subspace of $\,C[D]\,$.
\mp
In view of Lemma 2 it is clear that 
for $\,D_1,D_2\in{\cal D}\,$ the two spaces $\,C[D_1]\,$ and $\,C[D_2]\,$
are homeomorphic  if and only if there 
is a monotonic bijection from $\,[0,1]\,$ into $\,[0,1]\,$
which maps $\,D_1\,$ onto $\,D_2\,$.
Since, naturally, $\;|{\cal D}|=2^{\aleph_1}\;$
and since there are only $\,\bc\,$ monotonic mappings
from $\,[0,1]\,$ onto $\,[0,1]\,$, we can track down 
$\,2^{\aleph_1}\,$ mutually non-homeomorphic spaces 
in the family $\;\{\,C[D]\;|\;D\in{\cal D}\,\}\;$  
{\it provided that} $\,2^{\aleph_1}>\bc\,$.
Therefore, since both hypotheses 
$\,2^{\aleph_1}>\bc\,$ and $\,2^{\aleph_1}=\bc\,$ are consistent 
with ZFC, the only chance to prove Theorem 5 is to work 
in a model of ZFC where $\,2^{\aleph_1}=\bc\,$ holds
and hence  the Continuum Hypothesis fails.
(From $\,2^{\aleph_1}=\bc\,$ we derive $\,\aleph_1<\bc\,$
because $\,2^\bc>\bc\,$.)
Concerning the case $\,\kappa=\aleph_1\,$ in Theorem 4, 
we already have established 
that if $\,2^{\aleph_1}>\bc\,$ then
{\it there exist $\,2^{\aleph_1}\,$ mutually non-homeomorphic 
first countable linearly ordered continua of weight $\,\aleph_1\,$.}
However, if $\,2^{\aleph_1}=\bc\,$ then it is impossible 
to establish this enumeration result within the realm
of the special continua $\,C[D]\,$.
\mp
{\bf Proposition 3.} {\it It is consistent with {\rm ZFC}
that $\,2^{\aleph_1}=\bc\,$ and 
the family $\;\{\,C[D]\;|\;D\in{\cal D}\,\}\;$  contains
precisely three spaces up to homeomorphism.}
\mp
{\it Proof.} In [1] Baumgartner constructs a model 
of ZFC where $\,2^{\aleph_1}=\bc\,$
and all $\aleph_1$-dense subsets of $\,\R\,$ without extrema
are order-isomorphic. 
Thus if we define a partition $\,\{\,{\cal Y}_0,{\cal Y}_1,{\cal Y}_2\,\}\,$
of our family $\,{\cal D}\,$ via
\sp
\cl{$\;{\cal Y}_j\;:=\;\{\,D\in{\cal D}\;\,|\,\;|D\cap\{0,1\}|=j\,\}
\;\;(\,j\in\{0,1,2\}\,)$}
\sp
then within Baumgartner's universe for $\,j\in\{0,1,2\}\,$
and $\,D_1,D_2\in{\cal Y}_j\,$ there is always 
a monotonic bijection from $\,D_1\,$ onto $\,D_2\,$
which, of course, 
can be extended to one from $\,[0,1]\,$ onto $\,[0,1]\,$,
whence 
the spaces $\,C[D_1]\,$ and $\,C[D_2]\,$ are homeomorphic.
\sp
If $\;(D_1,D_2,D_3)\,\in\,
{\cal Y}_1\times{\cal Y}_2\times{\cal Y}_3\;$
then the three spaces $\,C[D_1]\,$ and $\,C[D_2]\,$ and $\,C[D_3]\,$
are  mutually non-homeomorphic because for $\,j\in\{0,1,2\}\,$
the total number of all endpoints $\,e\,$ in $\,C[D_j]\,$ 
where $\,\{e\}\,$ is a path component of $\,C[D_j]\,$
equals $\,2-j\,$, {\it q.e.d.}
\mp\sp
Since $\,2^{\aleph_1}=\bc\,$ implies $\,\aleph_1<\bc\,$,
Theorem 5 is an immediate consequence of 
Proposition 3 and the following proposition.
\mp
{\bf Proposition 4.} {\it If $\,\aleph_1<\bc\,$ and
$\,X\,$ is a first countable 
linearly ordered continuum of weight $\,\aleph_1\,$
such that the union of all singleton path components of $\,X\,$
is a separable subspace of $\,X\,$ then $\,X\,$ is homeomorphic
to the space $\,C[D]\,$ for some $\,D\in{\cal D}\,$.}
\mp
{\it Proof.}  Assume that $\,\aleph_1<\bc\,$ and
let $\,X\,$ be linearly ordered set which is  
a first countable linearly ordered continuum
with $\,w(X)=\aleph_1\,$. Let $\,{\cal P}_0\,$ 
denote the family of all non-singleton 
path components of $\,X\,$ and let
$\,{\cal P}_1\,$ denote the family of all singleton 
path components of $\,X\,$ 
and put $\;{\cal P}\,=\,{\cal P}_0\cup{\cal P}_1\;$ 
and assume that 
$\;\bigcup{\cal P}_1\,$ is a separable subspace of $\,X\,$.
Since $\,X\,$ is not second countable, $\,X\,$ is not homeomorphic to
$\,[0,1]\,$ and hence $\,X\,$ is not pathwise connected
and hence $\,|{\cal P}|\geq 2\,$.
\sp
The family $\,{\cal P}\,$ 
carries a natural ordering referring to the linear ordering of $\,X\,$.
Simply declare $\,A\in{\cal P}\,$ smaller than $\,B\in{\cal P}\,$
when every point in the set $\,A\,$ is smaller than every point
in the set $\,B\,$ within the linearly ordered set $\,X\,$.
(Referring to a standard terminology, see [11],
$\,{\cal P}\,$ is the {\it condensation} of $\,X\,$ with respect
to the equivalence relation $\,\sim\,$ where $\,x\sim y\,$
for distinct $\,x,y\in X\,$ when $\,w(I[x,y])=\aleph_0\,$.) 
By Lemma 5 all members of $\,{\cal P}\,$ are compact 
subspaces of $\,X\,$. (Hence each member of $\,{\cal P}_0\,$ 
is an order-isomorphic copy 
of the naturally or backwards ordered set $\,[0,1]\,$.)
Consequently, the linearly ordered set $\,{\cal P}\,$
is {\it densely ordered}. (Otherwise $\,X\,$ would not 
be connected.) If $\;\emptyset\not={\cal A}\subset {\cal P}\,$
then the supremum resp.~infimum of $\,{\cal A}\,$ 
is the unique path component which contains
the supremum resp.~infimum of the set $\,\bigcup{\cal A}\,$ in $\,X\,$.
Consequently, equipped with the order topology,
$\,{\cal P}\,$ is a linearly ordered continuum.
In particular, $\,|{\cal P}|\geq \bc\,$ and hence
$\,|{\cal P}|=\bc\,$ since the first countable continuum $\,X\,$ has 
size $\,\bc\,$. (Besides, $\,{\cal P}\,$ is first countable as well.)
We have $\,|{\cal P}_0|\leq\aleph_1\,$
because $\,w(X)\geq |{\cal P}_0|\,$ by Lemma 3.
Therefore, $\,|{\cal P}_1|=\bc\,$ since we assume $\,\aleph_1<\bc\,$.
\sp
For distinct $\,A,B\in{\cal P}\,$ consider the interval 
$\,I[A,B]\,$ in the linearly ordered set $\,{\cal P}\,$.
Naturally, the set $\,\bigcup I[A,B]\,$
is an interval in $\,X\,$. More precisely, if $\,A\,$ is smaller than 
$\,B\,$ in the linearly ordered set $\,{\cal P}\,$
then  $\,\bigcup I[A,B]\,$ equals the interval
$\,I[\min A,\max B]\,$ in $\,X\,$.
Since $\,|{\cal P}|=\bc\,$ and $\,{\cal P}\,$ is a linearly 
ordered continuum, we always have $\,|I[A,B]|=\bc\,$.
Furthermore we claim that
\sp
(6.2)\quad $\;|{\cal P}_0\cap I[A,B]|\,=\,\aleph_1\;$
{\it whenever $\;A,B\in{\cal P}\;$ and $\;A\not=B\,$.}
\sp
Assume indirectly that for distinct $\;A,B\in{\cal P}\;$
the set $\;{\cal P}_0\cap I[A,B]\;$ is only countable.
Then $\,J=\bigcup I[A,B]\,$ is an interval in $\,X\,$
which is partitioned into countably many copies of $\,[0,1]\,$
and the one separable space $\;J\cap \bigcup{\cal P}_1\,$.
Therefore, $\,J\,$ is a separable linearly ordered 
continuum and hence $\,J\,$ is homeomorphic to $\,[0,1]\,$ by (6.1).
But then $\,J\subset C\,$ for some path component $\,C\,$ of $\,X\,$
in contradiction to $\,A,B\subset J\,$ and
$\,A,B\in{\cal P}\,$ and $\,A\not=B\,$.
\sp
From (6.2) and $\,|I[A,B]|=\bc\,$ for distinct $\,A,B\in{\cal P}\,$ 
and $\,\aleph_1<\bc\,$ we conclude that 
\eject
\sp
(6.3)\quad $\;|{\cal P}_1\cap I[A,B]|\,=\,\bc\;$
{\it whenever $\;A,B\in{\cal P}_0\;$ and $\;A\not=B\,$.}
\sp
As a consequence of (6.3), the order topology of
$\,{\cal P}_1\,$ coincides with the subspace topology
of the continuum $\,{\cal P}\,$
and $\,{\cal P}_1\,$ is a dense point set in the space $\,{\cal P}\,$.
Clearly, the mapping $\;\{x\}\mapsto x\;$ is an order-isomorphism 
from the linearly ordered 
set $\,{\cal P}_1\,$ onto the subset $\,\bigcup{\cal P}_1\,$ of $\,X\,$
which is a separable subspace of $\,X\,$.
Consequently, $\,{\cal P}_1\,$ is a dense and separable
subspace of $\,{\cal P}\,$ and hence the space $\,{\cal P}\,$ 
itself must be separable.
Since $\,{\cal P}\,$ is a separable linearly ordered continuum,
by (6.1) and Lemma 2 there is an order-isomorphism $\,\varphi\,$
from $\,{\cal P}\,$ onto the (naturally or backwards 
ordered) set $\,[0,1]\,$. And $\,\varphi\,$  must map $\,{\cal P}_0\,$
onto a set $\,D\,$ in the family $\,{\cal D}\,$ in view of (6.2).
Since each member of $\,{\cal P}_0\,$ is an order-isomorphic copy
of $\,[0,1]\,$, in an obvious and natural way
the order-isomorphism $\,\varphi\,$ induces a monotonic 
bijection from the linearly ordered set $\,X\,$ 
onto the lexicographically ordered set $\,C[D]\,$.
Thus by Lemma~2 the linearly 
ordered continuum $\,X\,$ and $\,C[D]\,$ are homeomorphic, {\it q.e.d.}
\bp
{\bff 7. Proof of Theorem 4}
\mp
Throughout this section 
let us call a linearly ordered continuum $\,L\,$  {\it short} when 
every separable subspace is second countable.
In other words a linearly ordered continuum is short
if and only if {\it the closure of any countable point set is 
a second countable subspace.}
Clearly if $\,L\,$ is short then $\,I[x,y]\,$ is short 
whenever $\,x,y\in L\,$ and $\,x\not=y\,$.
\mp
{\bf Lemma 6.} {\it There exists a short, first countable 
linearly ordered continuum $\,C\,$ such that $\,\hat w(C)=w(C)=\aleph_1\,$.}
\mp
{\it Proof.} Let $\,A\,$ be a {\it densely ordered Aronszajn line}
and let $\,B\,$ be a compact Dedekind-completion of $\,A\,$.
By virtue of Proposition 3.6 in [8] Ch.~6, \S 3,
$\,B\,$ is a short and first countable linearly ordered continuum
and $\,w(B)\geq \aleph_1\,$.
Consequently, $\,w(B)=\aleph_1\,$ since $\,|A|=\aleph_1\,$
and $\,A\,$ is dense in $\,B\,$.
The path components of $\,B\,$ are compact by virtue of Lemma 5.
Therefore, 
the canonically ordered family $\,C\,$ of all
path components of $\,B\,$ is 
a first countable 
linearly ordered continuum with $\,\hat w(C)\geq \aleph_1\,$ and hence 
$\,\hat w(C)=w(C)=\aleph_1\,$. 
Finally, the continuum $\,C\,$ is short. Indeed, 
let $\,S\subset C\,$ be countable. Choose for every 
path component $\,P\,$ of $\,B\,$ with $\,P\in S\,$
a nonempty countable set $\,T_P\subset P\,$ 
dense in $\,P\,$. (If $\,P\,$ is a singleton 
then $\,T_P=P\,$. If $\,P\,$ is  not a singleton then 
$\,P\,$ is a copy of $\,[0,1]\,$ by Lemma 5 and 
$\,T_P\,$ can be chosen countably infinite.)
Let $\,T\,$ be the closure of the countable set 
$\;\bigcup_{P\in S}T_P\;$ in $\,B\,$.
Obviously, $\;{\cal P}\,=\,\{\,P\in C\;|\;P\cap T\not=\emptyset\,\}\;$
is the closure of $\,S\,$ in $\,C\,$. And $\,{\cal P}\,$
is second countable in $\,C\,$ since 
$\,T\,$ is second countable in $\,B\,$, {\it q.e.d.}   
\mp
{\it Remark.} If we start the previous proof 
with a densely ordered Aronszajn line $\,A\,$ such that 
$\,\hat w(A)=\aleph_1\,$ then
the condensation which creates $\,C\,$ from $\,B\,$
is superfluous because $\,\hat w(B)=\aleph_1\,$ follows from 
$\,\hat w(A)=\aleph_1\,$ and $\,C=B\,$ fits already.
However, while it is trivial to turn an Aronszajn line $\,A\,$ 
into a {\it densely ordered} Aronszajn line 
(simply consider the lexicographically ordered set $\,A\times \Q\,$),
some explanations have to be carried out in order to define a 
densely ordered Aronszajn line $\,A\,$ with  
$\,\hat w(A)=\aleph_1\,$ correctly.
\mp\sp 
In the following let $\,C\,$ be a continuum as provided by Lemma 6.
Furthermore let $\,K\,$ be a 
first countable linearly ordered continuum with $\,\hat w(K)=w(K)=\aleph_1\,$
and no compact interval in $\,K\,$ is short.
Such a continuum $\,K\,$ can easily be constructed.
Indeed, referring to Section 4, 
choose any dense $\,D\subset [0,1]\,$ such that $\,1\in D\,$
and $\,|D|=\aleph_1\,$ and put $\,K=L[D]\,$ for $\,L=[0,1]\,$.
Consider the countable set 
$\,A\,$ of all points $\,(r,0,0,0,...\,)\,$ with $\,r\in \Q\cap[0,1]\,$.
Obviously, the closure $\,B\,$ of $\,A\,$ in $\,L[D]\,$ 
consists of the points $\,(x,0,0,0,...\,)\;$
and $\,(y,1,1,1,...\,)\,$ with $\,0\leq x\leq 1\,$
and $\,1\not=y\in D\,$. Hence $\,w(B)\geq |D|=\aleph_1\,$ by Lemma~4.
Consequently, since every interval in $\,K\,$
contains copies of $\,K\,$, if $\,x,y\in K\,$ are distinct
then $\,I[x,y]\,$ cannot be short.
\mp
Let us call two linearly ordered continua $\,X_1,X_2\,$ 
{\it extremely incomparable} if and only if for $\,\{i,j\}=\{1,2\}\,$
no interval in $\,X_i\,$ is homeomorphic to a subspace 
of $\,X_j\,$. In the following it is essential that
$\,C\,$ and $\,K\,$ are extremely incomparable.
\mp
Now we are ready to prove Theorem 4. We assume $\;\aleph_0<\kappa\leq \bc\;$
and distinguish the cases $\,2^\kappa>\bc\,$ and $\,2^\kappa=\bc\,$.
Assume firstly that $\,2^\kappa>\bc\,$.
Let $\,{\cal D}\,$ be a family of subsets of $\;[1,2]\setminus\Q\;$ 
which are dense in $\,[1,2]\,$ and of size $\,\kappa\,$ 
such that $\,|{\cal D}|=2^\kappa\,$ 
and for distinct $\,D_1,D_2\,$
there is never a monotonic bijection $\,f\,$ from $\,[1,2]\,$
onto $\,[1,2]\,$ with $\,f(D_1)=D_2\,$.
(Such a family $\,{\cal D}\,$ exists under the assumption 
$\,2^\kappa>\bc\,$ since there are $\,2^\kappa\,$ 
dense subsets of $\;[1,2]\setminus\Q\;$ but
only $\,\bc\,$ monotonic permutations
of the set $\,[1,2]\,$.)
Put $\;A\,:=\,([1,2]\cap\Q)\cup\{0,3\}\;$ 
and $\;B\,:=\,(]0,1[\cup]2,3[)\cap\Q\,$. 
For every $\,D\in{\cal D}\,$ consider
\mp
\cl{$X[D]\;:=\;([0,3]\setminus(A\cup B\cup D))\!\times\!\{0\}\,\cup\,
A\!\times\! C\,\cup\,(B\cup D)\!\times\! K$}
\mp
equipped with the lexicographic ordering.
So $\,X[D]\,$ is constructed from 
$\,[0,3]\,$  by replacing each point in $\,A\,$ with a copy of $\,C\,$
and each point in $\,B\cup D\,$ with a copy of $\,K\,$.
(Clearly, if $\;(x_1,x_2),(y_1,y_2)\in X[D]\;$ and $\,x_1\not= x_2\,$ 
then it is irrelevant whether $\,x_2\,$ and $\,y_2\,$ 
are comparable.) 
In view of Lemma 1 it is clear 
that the linearly ordered space $\,X[D]\,$ is compact and connected
and first countable. Since $\,A\cup B\cup D\,$ is dense in $\,[0,3]\,$
and $\,\hat w(K)=\aleph_1\,$ and $\,\hat w(C)=\aleph_1\,$,
we have $\;\hat w(X[D])=\aleph_1\,$. Since $\,w([0,3])=\aleph_0\,$
and $\,w(K)=w(C)=\aleph_1\,$ and $\,|A\cup B\cup D|=|D|=\kappa\,$, 
we have $\;w(X[D])=\kappa\,$.
\mp
We claim that the $\,2^\kappa\,$ spaces 
$\;X[D]\;(D\in{\cal D})\;$ are incomparable. 
Let $\,D_1,D_2\in{\cal D}\,$ and $\,f\,$ 
be a homeomorphism from $\,X[D_1]\,$ onto 
a subspace of $\,X[D_2]\,$.    
Then $\,f(X[D_1])\,$ is an interval in $\,X[D_2]\,$. 
Clearly, the subspaces $\,\{x\}\times C\,$ resp.~$\,\{y\}\times K\,$
of $\,X[D_i]\,$ are intervals homeomorphic to $\,C\,$ resp.~$\,K\,$.
Therfore, since $\,B\cup D\,$ is dense in $\,[0,3]\,$ and $\,C,K\,$ are 
extremely incomparable, if $\;a\in\{0,3\}\;$ then 
$\;f(\{a\}\times C)\,\subset\,\{a'\}\times C\;$ for some $\,a'\in A\,$
and if $\;a\,\in\,A\setminus\{0,3\}\;$ then
$\;f(\{a\}\times C)\,=\,\{a'\}\times C\;$ for some $\,a'\in A\,$.
Moreover, since
$\;A\cap\,]1,2[\;$ is order-isomorphic with $\,\Q\,$
and $\;A\cap[0,1]=\{0,1\}\;$ and $\;A\cap[2,3]=\{2,3\}\,$,
we have $\;\{f(\{1\}\times C),f(\{2\}\times C)\}\,=\,
\{\{1\}\times C,\,\{2\}\times C\}\,.$
As a consequence, 
$\,f\,$ induces a monotonic bijection $\,g\,$
from $\,[1,2]\,$ onto $\,[1,2]\,$ such that $\,g(D_1)=D_2\,$
and hence we obtain $\,D_1=D_2\,$. 
\sp
This proves the claim
and settles the case $\,2^\kappa>\bc\,$ in Theorem 4.
In order to settle the case $\,2^\kappa=\bc\,$
we will apply the following lemma.
\mp
{\bf Lemma 7.} {\it There exists a family $\,{\cal H}\,$
of countable subsets of $\,[0,1]\,$ 
such that $\,|{\cal H}|=\bc\,$ and for distinct 
$\,H_1,H_2\in {\cal H}\,$ and for $\,0\leq u<v\leq 1\,$
there is never a monotonic bijection from 
$\,H_1\,$ onto $\,H_2\cap[u,v]\,$.}
\sp
{\it Proof.} Similarly as in the proof of Proposition 1 
let $\,G\,$ be the family of 
all functions from $\,\N\,$ to $\,\{3,5\}\,$
and let $\,\zeta\,$ be the order-type of the naturally ordered set $\,\Z\,$.
Then choose for every $\,g\in G\,$
a set $\,H_g\subset[0,1[\;$ of order-type 
$\;7+\zeta+g(1)+\zeta+g(2)+\zeta+g(3)+\,\cdots\;$
and put $\;{\cal H}\,=\,\{\,H_g\;|\;g\in G\,\}\,$, {\it q.e.d.}
\mp
Now in order to conclude the proof of Theorem 4 
assume that $\;\aleph_0<\kappa\;$ and $\,2^\kappa=\bc\,$.
Let $\,{\cal H}\,$ be a family as provided by Lemma 7.
Choose for every (countable) $\,H\in{\cal H}\,$ 
a dense set $\,D(H)\subset [0,1]\,$ 
with $\,H\cap D(H)=\emptyset\,$ and $\,|D(H)|=\kappa\,$
and consider 
\mp
\cl{$Y[H]\;:=\;([0,1]\setminus(H\cup D[H])\!\times\!\{0\}\,\cup\,
H\!\times\! C\,\cup\,D[H]\!\times\! K$}
\mp
equipped with the lexicographic ordering.
Similarly as above, $\,Y[H]\,$ is 
a first countable linearly ordered continuum 
and $\;\hat w(Y[H])=\aleph_1\;$
and $\;w(Y[H])=\kappa\,$.
The $\,\bc\,$ spaces
$\;Y[H]\;(H\in{\cal H})\;$ are incomparable 
because, by similar arguments as above, for $\,H_1,H_2\in {\cal H}\,$
any homeomorphism between $\,Y[H_1]\,$ and 
a subcontinuum of $\,Y[H_2]\,$ 
induces a monotonic bijection from $\,H_1\,$ onto 
$\,H_2\cap[u,v]\,$ for some $\,(u,v)\,$ with $\,0\leq u<v\leq 1\,$
and this is only possible if $\,H_1=H_2\,$.
This concludes the proof of Theorem 4.

\bp
{\bff 8. Proof of Theorem 7}
\mp
For abbreviation put $\,\exp\kappa\,:=\,2^\kappa\,$.
As usual, if $\,X\,$ is a Hausdorff space then
$\,c(X)\,$ denotes the {\it cellularity} of $\,X\,$.
(So $\,c(X)\,$ is the supremum of the sizes $\,|{\cal U}|\,$
of all families $\,{\cal U}\,$ of mutually disjoint open 
subsets of $\,X\,$.)
Trivially, $\,c(X)\leq d(X)\leq w(X)\,$.
Although $\;|X|\leq \exp\max\{c(X),\chi(X)\}\;$ for any 
Hausdorff space $\,X\,$ (see [4] 2.15),
in general it is not possible to estimate the size or the weight
or the density 
of a space with an upper bound depending on the cellularity only.
In fact, the gap 
between $\,c(X)\,$ and $\,|X|,w(X),d(X)\,$ 
may be arbitrarily large 
because for the Hilbert cube $\,[0,1]^\kappa\,$ 
we always have 
$\;d([0,1]^\kappa)=\log\kappa\;$
and $\;w([0,1]^\kappa)=\kappa\,$
and $\;|[0,1]^\kappa|=2^\kappa\,$
{\it but} $\;c([0,1]^\kappa)=\aleph_0\,$, see
[4] 5.3, 5.5 and [2] 2.7.10.d.
(In comparison, for the linearly ordered continuum $\,L=[0,1]^{[0,\kappa[}\,$
we also have
$\;|L|=2^\kappa\,$, but we have
$\;d(L)=w(L)=\bc^{<\kappa}\,$.
This is the main reason why it is more difficult,
even without considering incomparability,  
to count {\it linearly ordered} continua up to homeomorphism,
than to accomplish the enumeration results in [9] 
concerning continua {\it only}.)
However, the following two lemmas show
that there are two special cases where an estimate 
depending on the cellularity only is indeed possible.
(Note that $\;|X|\leq \exp w(X)\;$
and $\;|X|\leq \exp\exp d(X)\;$ for every Hausdorff space $\,X\,$, see 
[4] 2.2 and 2.4.)
\mp
{\bf Lemma 9.} {\it If $\,X\,$ is a 
compact and hereditarily normal space then 
$\,w(X)\leq \exp c(X)\,$.} 
\mp
{\bf Lemma 10.} {\it If $\,X\,$ is a linearly ordered space
then either $\,d(X)=c(X)\,$ or $\,d(X)=c(X)^+\,$.}
\mp
Lemma 9 is an immediate consequence
of the elementary estimate [4] 2.7 and \v Sapirovskii's deep result
[4] 3.21. Lemma 7 is a deep result due to Hajnal and Juh sz (cf.~[2] 3.12.4).
\mp
{\bf Lemma 11.} {\it If $\,\lambda\,$ is not a regular limit and 
$\,X\,$ is a space 
with $\,c(X)=\lambda\,$ 
then $\,X\,$ contains $\,\lambda\,$ disjoint open sets.}
\mp
Lemma 11 is trivially true
if $\,\lambda\,$ is a successor cardinal.
A proof of Lemma 11 for singular $\,\lambda\,$ is given in [4] 4.1. 
(Lemma 11 would be false for regular limits $\,\lambda\,$,
see [4] 7.6.) 
\mp
The second statement of Theorem 7 is an immediate consequence 
of the following theorem.
\mp
{\bf Theorem 11.} {\it If $\,\kappa\,$ is singular and 
$\,\lambda\leq \kappa\,$
and $\,K\,$ is a linearly ordered continuum such that $\,\hat w(K)=\kappa\,$
and no point in $\,K\,$ has a local basis of size smaller than $\,\lambda\,$
then $\,|K|\geq\kappa^\lambda\,$.}
\mp
{\it Proof.} Let $\,{\cal J}\,$ be the family of all
intervals in $\,K\,$. (So every $\,J\in{\cal J}\,$ is an {\it infinite}
compact and connected subset of $\,K\,$.) 
If $\,\emptyset\not={\cal G}\subset{\cal J}\,$
and $\,{\cal G}\,$ is a {\it chain}
(i.e.~$\,I\subset J\,$ or $\,J\subset I\,$ whenever $\,I,J\in{\cal G}\,$)
then $\;\bigcap{\cal G}\not=\emptyset\;$ since $\,K\,$ is compact.
If additionally $\,|{\cal G}|<\lambda\,$ then
the intersection of the intervals in $\,{\cal G}\,$
is not only nonempty but an interval itself. Because otherwise
$\,\bigcap{\cal G}\,$ would be a singleton $\,\{x\}\,$
and the family $\,{\cal B}\,$ of the interiors of the intervals
$\,J\in{\cal G}\,$ would be a local basis at $\,x\,$ 
with $\,|{\cal B}|<\lambda\,$.
(Here it is essential that $\,K\,$ is a linearly ordered
space. An analogous argument would not work for arbitrary continua!)
\sp
Consider any linearly ordered continuum $\,J\in{\cal J}\,$.
Then $\,d(J)=w(J)=\hat w(K)=\kappa\,$.
Since $\,\kappa\,$ is a limit cardinal,
we have $\,c(J)=d(J)\,$ by Lemma 10.
Consequently, since $\,\kappa\,$ is singular,
by virtue of Lemma 11 the interval $\,J\,$ 
contains $\,\kappa\,$ mutually disjoint intervals.
Therefore, we can define for
each $\,J\in{\cal J}\,$ a function $\,f[J]\,$ from $\,[0,\kappa[\;$
into $\,{\cal J}\,$ such that $\,f[J](\alpha)\subset J\,$ 
for every $\,\alpha\in[0,\kappa[\;$ and 
$\;f[J](\alpha)\cap f[J](\beta)=\emptyset\;$
whenever $\,\alpha<\beta<\kappa\,$.
\sp
Let $\;(\alpha_\beta)_{\beta<\lambda}\;$ be a $\lambda$-sequence 
of ordinals $\;\alpha_\beta\in[0,\kappa[\,$.
In the following, by transfinite induction 
we define for every ordinal $\,\gamma<\lambda\,$
an interval 
$\;I((\alpha_\beta)_{\beta\leq\gamma})\;$ 
such that $\,I((\alpha_\beta)_{\beta\leq\delta})\,$
is a proper subset of $\,I((\alpha_\beta)_{\beta\leq\gamma})\,$
whenever $\;\gamma<\delta<\lambda\,$.
\sp
Assume that for $\,\gamma<\lambda\,$
a descending chain $\;\{\,I((\alpha_\beta)_{\beta\leq\delta})\;|
\;\delta<\gamma\,\}\;$ of intervals has already been defined.
Then put $\,\Delta_0=K\,$ and if $\,\gamma>0\,$ then put
$\;\Delta_\gamma\,:=\,
\bigcap_{\delta<\gamma}I((\alpha_\beta)_{\beta\leq\delta})\,$.
The set $\,\Delta_\gamma\,$ is 
an interval since $\;|[0,\gamma[|<\lambda\,$.
(If $\,\gamma=\eta+1\,$ for some ordinal $\,\eta\,$
then, trivially, 
$\;\Delta_\gamma=I((\alpha_\beta)_{\beta\leq\eta})\,$.)
This allows us to proceed and 
define $\;I((\alpha_\beta)_{\beta\leq\gamma})\,:=\,
f[\Delta_\gamma](\alpha_{\gamma})\,$.
(Then, obviously, the family 
$\;\{\,I((\alpha_\beta)_{\beta\leq \delta})\;
|\;\delta\leq\gamma\,\}\;$ is a descending chain of intervals.)
\sp
By construction, for every   $\lambda$-sequence 
$\;a=(\alpha_\beta)_{\beta<\lambda}\;$ 
of ordinals $\;\alpha_\beta\in[0,\kappa[\;$
the family $\;{\cal G}(a)\,:=
\{\,I((\alpha_\beta)_{\beta\leq\gamma})\;|\;\gamma<\lambda\,\}\;$
is a chain of intervals and hence
$\;\bigcap{\cal G}(a)\not=\emptyset\;$ since $\,K\,$ is compact.
And by the definition of the chains $\,{\cal G}(a)\,$
it is clear that $\,\bigcap{\cal G}(a_1)\,$ and $\,\bigcap{\cal G}(a_2)\,$ 
are always disjoint for distinct $\lambda$-sequences $\,a_1\,$ and $\,a_2\,$.
Consequently, $\,|K|\,$ is not smaller than the total amount of all
$\lambda$-sequences of ordinals smaller than $\,\kappa\,$.
In other words, $\;|K|\geq \kappa^\lambda\,$, {\it q.e.d.}
\mp
{\it Remark.} As a trivial but noteworthy  consequence of Theorem 11, 
if $\,K\,$ is a balanced linearly 
ordered continuum such that $\,|K|<\exp w(K)\,$
then the character of some point in $\,K\,$ must be smaller than $\,w(K)\,$.
In our proof of Theorem 6 (and Theorem 3) it was essential 
that not all points have the same character.
\mp
The arguments in the previous proof can be used to 
derive a noteworthy theorem about locally compact Hausdorff spaces
which also implies 
the second statement of Theorem 7
and which we will apply in order to prove the first 
statement of Theorem 7.
Note that, by virtue of Lemma 10, for a limit cardinal $\,\kappa\,$
and a linearly ordered space $\,Y\,$ we have $\,d(Y)\geq \kappa\,$
if and only if $\,c(Y)\geq \kappa\,$.
\mp
{\bf Theorem 12.} {\it If $\,\kappa\,$ is not a regular limit 
and $\,X\,$ is a locally compact Hausdorff space 
such that $\;c(U)\geq\kappa\;$ for each 
infinite open subspace $\,U\,$ of $\,X\,$ 
then $\;|X|\geq \kappa^{\aleph_0}\,$.}
\mp
{\it Proof.} The space $\,X\,$ contains only finitely many 
isolated points because
if $\,A\subset X\,$ and $\,|A|=\aleph_0\,$
and all points in $\,A\,$ are isolated in $\,X\,$ then
$\,A\,$ is open in $\,X\,$ and 
$\,c(A)=|A|=\aleph_0\,$ contrarily to $\,c(A)\geq\kappa>\aleph_0\,$.
So we may assume that no point in $\,X\,$ is isolated. 
For simply speaking 
let us call $\,G\,$ {\it good} when $\,G\,$ 
is a compact subset of $\,X\,$ which is the closure of some nonempty open set. 
By definition, any good set contains a nonempty open set.
Since $\,X\,$ is locally compact, any nonempty open set contains a good set.
Therefore, from the assumption $\;c(U)\geq\kappa\;$ for each 
infinite open subspace $\,U\,$ of $\,X\,$ 
(and by applying Lemma~11 in case that $\;c(U)=\kappa\,$)
we conclude that {\it any good set contains $\,\kappa\,$
mutually disjoint good sets.}
This enables us
to choose for each $\,n<\omega\,$ 
and every $(n+1)$-tuple $\,(\alpha_0,...,\alpha_n)\,$ of ordinals 
$\,\alpha_i\,$ smaller than $\,\kappa\,$
a good set 
$\;G(\alpha_0,...,\alpha_n)\;$ 
such that always $\;G(\alpha_0,...,\alpha_n)\cap 
G(\beta_0,...,\beta_n)\not=\emptyset\;$  for distinct 
$\,(\alpha_0,...,\alpha_n)\,$ and $\,(\beta_0,...,\beta_n)\,$ and
\sp 
\cl{$G(\alpha_0)\,\supset\, G(\alpha_0,\alpha_1)\,\supset\,
G(\alpha_0,\alpha_1,\alpha_2)\,\supset\,\cdots$}
\sp
for every sequence $\;(\alpha_0,\alpha_1,\alpha_2,...\,)\,.$
For every sequence $\;a=(\alpha_0,\alpha_1,\alpha_2,...\,)\;$
the set $\;S(a)\,=\,\bigcap_{n<\omega} G(\alpha_0,...,\alpha_n)\;$
is not empty since $\,G(\alpha_0)\,$ is compact.
And clearly for distinct sequences $\;a_1,a_2\;$ 
the sets $\,S(a_1)\,$ and $\,S(a_2)\,$ are disjoint.
Therefore $\,|X|\,$ is not smaller than the total amount
of all these sequences. In other words, $\;|X|\geq \kappa^{\aleph_0}\,$,
{\it q.e.d.}
\mp\sp
Finally, the first statement of Theorem 7 is an immediate consequence 
of the following interesting and noteworthy theorem.
\mp
{\bf Theorem 13.} {\it If $\,\kappa\,$ is a strong limit with 
$\,\kappa^{\aleph_0}>\kappa\,$
and $\,X\,$ is a locally compact 
and hereditarily normal space of size  $\,\kappa\,$
then $\;|U|<\kappa\;$ for some nonempty open set $\,U\subset X\,$.}
\mp
{\it Proof.} 
Since there is nothing to show if 
$\,X\,$ contains finite nonempty open sets, we may assume that
all nonempty open subsets of $\,X\,$ are infinite.
Then $\,\kappa>\aleph_0\,$ because
any countable locally compact Hausdorff space
must have isolated points. In particular,
the strong limit $\,\kappa\,$  of countable cofinality
is a singular cardinal.
Since $\,\kappa^{\aleph_0}>\kappa\,$,
by applying Theorem 12 there
exists an infinite open set $\,V\subset X\,$ with 
$\,c(V)<\kappa\,$. Trivially, for each infinite 
open set $\,W\subset V\,$ we have $\,c(W)<\kappa\,$ as well.
Hence we can find an infinite open set $\,U\subset X\,$ 
such that the closure $\,\overline U\,$ of $\,U\,$ in $\,X\,$
is compact and $\,c(U)<\kappa\,$.
Since $\,U\,$ is dense in $\,\overline U\,$, we have $\,c(\overline U)=c(U)\,$
(cf.~[4] 2.6). Since $\,\kappa\,$ is a strong limit, 
from $\,c(\overline U)<\kappa\,$ we derive
$\;\exp\exp c(\overline U)<\kappa\,$. 
Since $\,{\overline U}\,$ is compact and hereditarily normal,
via the basic estimate 
$\,|\overline U|\leq \exp w(\overline U)\,$ and Lemma 9 we derive
$\;|\overline U|\leq \exp\exp c(\overline U)\,$.
Consequently, $\;|U|\leq |{\overline U}|<\kappa\,$, {\it q.e.d.}
\bp
{\it Remark.} For a proof of Theorem 7 (where only
linearly ordered spaces are considered)
one can avoid Lemma 9 because for a linearly ordered
space $\,X\,$ we have  $\;\chi(X)\leq c(X)\;$ (cf.~[2] 3.12.4)
and hence  $\;w(X)\leq |X|\leq \exp c(X)\,$.
\bp
{\bff 9. Proof of Theorem 8}  
\mp
In the following we need 
{\it Cantor derivatives} and {\it signature sets}. 
If $\,X\,$ is a Hausdorff space and $\,\xi\,$
is an ordinal number and $\,A\subset X\,$ then 
$\,A^{(\xi)}\,$ is the $\xi$-{\it th derivative} of the point set $\,A\,$.
($\,A^{(0)}=A\,$ and $\,A^{(\alpha+1)}\,$
is the set of all limit points of $\,A^{(\alpha)}\,$ 
and if $\,\beta>0\,$ is a limit ordinal then
$\;A^{(\beta)}\,:=\,\bigcap\,\{\,A^{(\alpha)}\;|\;
0<\alpha<\beta\,\}\,$.)
It is essential that if $\,\omega^\xi\,$ denotes the
{\it ordinal power} with basis $\,\omega\,$ and exponent $\,\xi\,$
then in the compact linearly 
ordered space $\;X=[0,\omega^\xi]\;$ we have 
$\;[0,\omega^\xi[^{(\alpha)}=[0,\omega^\xi]^{(\alpha)}
=\{\omega^\xi\}\;$ if and only if $\,\alpha=\xi\,$.
(Note that $\;|[0,\omega^\xi]|=|[0,\xi]|\;$ whenever $\,\xi\geq\omega\,$.)
If $\,X\,$ is a Hausdorff space then 
$\;\Delta(X)\,=\,\bigcup\,\{\,A\subset X\;|\;A\subset A^{(1)}\,\}\;$ 
is {\it the maximal dense-in-itself point set}.
The set
\mp
\cl{$\Sigma(X)\;:=\;
\big\{\,\alpha\in [0,|X|^+]\;\,\big|\,\;
\big((X\setminus \Delta(X))^{(\alpha)}\setminus
(X\setminus\Delta(X))^{(\alpha+1)}\big)\cap \Delta(X)\,
\not=\,\emptyset\,\big\}$}
\mp
may be regarded as a signature set since 
two spaces $\,X_1,X_2\,$ cannot be homeomorphic
if $\;\Sigma(X_1)\not=\Sigma(X_2)\,$.
\mp\sp
The clue in order to settle Theorem 8 is the following interesting theorem.
\mp 
{\bf Theorem 14.} {\it Let $\,K\,$ be a linearly ordered
continuum. 
Then for every $\,\kappa\,$ there exists 
a family $\,{\cal Y}[K,\kappa]\,$ of mutually non-homeomorphic 
compact linearly ordered spaces such that 
$\;|{\cal Y}[K,\kappa]|=2^\kappa\;$ and if $\,Y\in{\cal Y}[K,\kappa]\,$
then $\,Y\,$ has precisely $\,\kappa\,$ components 
and all components of $\,Y\,$ are homeomorphic to $\,K\,$, whence
$\;w(Y)=\max\{\kappa, w(K)\}\;$ and $\,\hat w(Y)=\hat w(K)\,$
and $\;|Y|=\max\{\kappa,|K|\}\,$.}
\mp
It is worth mentioning that 
Theorem 14 has an interesting consequence 
about compact sets in the number line.
Any second countable linearly ordered 
space is homeo\-morphic to a subspace of $\,\R\,$
(cf.~[2] 4.2.9, 6.3.2.c). So if we put
$\,K=[0,1]\,$ and $\,\kappa=\aleph_0\,$ in Theorem~14 
then we immediately obtain the following enumeration result.
\mp
{\bf Corollary 2.} {\it  There are $\,\bc\,$ mutually non-homeomorphic 
compact subspaces of $\,{\Bbb R}\,$
without singleton components.}
\mp
{\it Remark.} By  Corollary 2, up to homeomorphism the number line $\,\R\,$ 
contains precisely $\,\bc\,$ {\it closures of open sets}.
On the other hand, $\,\R\,$ contains 
precisely $\,\aleph_0\,$ {\it open} subspaces up to homeomorphism.
(Because for every cardinal $\,\theta\leq\aleph_0\,$ 
there is essentially one and only one
open subspace of $\,\R\,$ with precisely $\,\theta\,$ components.)
This cardinal discrepancy 
is specific for the space $\,\R\,$
and vanishes in the spaces $\,\R^n\,$ for dimensions $\,n>1\,$
(cf.~[5] 8.1). There is another noteworthy cardinal discrepancy.
If $\,K\,$ is a compact subspace of $\,{\Bbb R}\,$
without singleton components then the boundary of $\,K\,$ in $\,\R\,$
is a compact and {\it countable} set. However, in view of [10]
there exist precisely $\,\aleph_1\,$ 
compact and countable Hausdorff spaces up to homeomorphism.
(And $\,\aleph_1<\bc\,$ is consistent with ZFC.)
Moreover, referring to the remark below the proof of Theorem 14, one can 
track down $\,\bc\,$ mutually non-homeomorphic compact subspaces of $\,{\Bbb R}\,$
without singleton components such that all their boundaries in $\,\R\,$
are homeomorphic subspaces of $\,\R\,$.
\bp
Now in order prove Theorem 14
let $\,K\,$ be a linearly ordered set 
such that, with respect to the order topology, 
$\,K\,$ is a continuum. For any ordinal number $\,\xi>0\,$
let $\,Z_\xi\,$ be a set linearly ordered by $\,\preceq\,$ such that 
for a unique point $\,z_\xi\in Z_\xi\,$
the set $\;\{\,x\in Z_\xi\;|\;x\preceq z_\xi\,\}\;$
is order-isomorphic to the well-ordered set $\,[0,\omega^\xi]\,$
and the set $\;\{\,x\in Z_\xi\;|\;x\succeq z_\xi\,\}\;$
is order-isomorphic to the anti-well-ordered set $\,\Z\setminus\N\,$.
\mp
Now consider the set 
$\;Z_\xi(K)\,:=\,Z_\xi\times K\;$
equipped with the lexicographic ordering.
So $\,Z_\xi(K)\,$ is constructed from $\,Z_\xi\,$
by replacing each point in $\,Z_\xi\,$ with a copy of $\,K\,$.
It goes without saying that $\,Z_\xi(K)\,$ is a compact 
linearly ordered space whose components are all homeomorphic with $\,K\,$.
Among these components the component $\,\{z_\xi\}\times K\,$
is labeled as the unique component $\,C\,$ of $\,Z_\xi(K)\,$
where both endpoints are limit points of the set $\,Z_\xi(K)\setminus C\,$.
Trivially, $\;|Z_\xi(K)|=\max\{|K|,|\xi|\}\,$.
Clearly, $\;w(Z_\xi(K))=\max\{w(K),|\xi|\}\,$.
\mp
Let $\,{\cal S}_\kappa\,$ denote the family of all
sets $\,S\subset [2,\kappa[\;$ such that 
$\,|S|=\kappa\,$, whence  $\,|{\cal S}_\kappa|=2^\kappa\,$.
For each $\,S\in{\cal S}_\kappa\,$ let 
$\,Y[S]\,$ be constructed from the well-ordered 
set $\,[0,\kappa]\,$ by replacing each point $\,\xi\in S\,$
with an isomorphic copy of the linearly ordered
set $\,Z_\xi(K)\,$ and by replacing each point $\,\xi\not\in S\,$
with an isomorphic copy of the linearly ordered
set $\,K\,$. So $\,Y[S]\,$ is a linearly ordered set 
with a minimum and a maximum and
every nonempty subset of $\,Y[S]\,$ has a supremum.
Hence $\,Y[S]\,$ is compact in the order-topology.
Obviously, all components of the compact linearly ordered space $\,Y[S]\,$
are homeomorphic to $\,K\,$.
Clearly, $\;|Y[S]|=\max\{\kappa,|K|\}\;$ 
and $\;w(Y[S])=\max\{\kappa,w(K)\}\,$.
So we put $\;{\cal Y}[K,\kappa]\,:=\,\{\,Y[S]\;|\;S\in{\cal S}_\kappa\,\}\;$
and are finished by showing that $\,Y[S_1]\,$ and $\,Y[S_2]\,$ are 
never homeomorphic for distinct 
$\,S_1,S_2\in{\cal S}_\kappa\,$.            
\mp
Let $\,{\cal C}[S]\,$ denote the family of all components
of the space $\,Y[S]\,$ and for each $\,C\in{\cal C}[S]\,$
let $\,E(C)\,$ be the set of the two endpoints of $\,C\,$.
(Topologically, $\,E(C)\,$ consists of the two noncut points
of the space $\,C\,$.)
Let $\,{\cal C}_2[S]\,$ denote the family of all $\,C\in{\cal C}[S]\,$
such that 
{\it both} endpoints of $\,C\,$ are limit points of the set 
$\;Y[S]\setminus C\,$. Obviously, the members
of $\,{\cal C}_2[S]\,$ are precisely the labeled 
components $\,\{z_\xi\}\times K\,$ with $\,\xi\in S\,$.
With respect to the linear ordering of $\,Y[S]\,$, 
for each $\;C\,\in\,{\cal C}[S]\setminus{\cal C}_2[S]\;$ 
either the left endpoint or no endpoint is a limit point 
of $\;Y[S]\setminus C\,$. 
For any labeled component $\;C\,=\,\{z_\xi\}\times K\;$
and the subspace $\;Z_\xi(K)\;$ of $\,Y[S]\,$ around $\,C\,$
the right endpoint of $\,C\,$ lies in 
$\,Z_\xi(K)^{(1)}\,$ but not in 
$\,Z_\xi(K)^{(2)}\,$, while (since $\,\xi\geq 2\,$)
the left endpoint of $\,C\,$ lies in $\,Z_\xi(K)^{(2)}\,$.
Therefore, if $\,C\in{\cal C}_2[S]\,$
and $\;R_C\,=\,\bigcup\,\{\,E(G)\;|\;C\not=G\in{\cal C}[S]\,\}\;$
then the right endpoint of $\,C\,$ lies in 
$\,R_C^{(1)}\setminus R_C^{(2)}\,$
and the left endpoint of $\,C\,$ lies in $\,R_C^{(2)}\,$.
Consequently, whenever $\,S\in {\cal S}_\kappa\,$, 
for the subspace 
\mp
\cl{$\;R[S]\;:=\;\bigcup\,\{\,E(C)\;|\;
C\,\in\,{\cal C}[S]\setminus {\cal C}_2[S]\,\}\;\cup\;
\bigcup{\cal C}_2[S]\;$}
\mp
of $\,Y[S]\,$ 
we have $\;\Sigma(R[S])\setminus\{1\}\,=\,S\;$
and this concludes the proof of Theorem 14.
\bp
{\it Remark.} Put $\,K=[0,1]\,$ and $\,\kappa=\aleph_0\,$.
Then every space $\,Y[S]\,$ is homeomorphic to a compact subspace $\,C_S\,$ 
of $\,\R\,$ without singleton components. The members of $\,{\cal E_\kappa}\,$ 
are the infinite sets $\;S\,\subset\,\N\setminus\{1\}\;$ 
and, of course, $\;|{\cal S}_\kappa|=\bc\,$.
The spaces $\;C_S\;(S\in{\cal S}_\kappa)\;$ are mutually non-homeomorphic.
However, if $\,B_S\,$ is the boundary of $\,C_S\,$ in $\,\R\,$,
whence $\,B_S\,$ is compact and countable, then all spaces 
$\;B_S\;(S\in{\cal S}_\kappa)\;$ are homeomorphic!
This is an immediate consequence of the famous classification 
theorem concerning compact and countable Hausdorff spaces in [10]
since $\,B_S^{(\omega)}\,$ is a singleton for every 
$\,S\in{\cal S}_\kappa\,$.
This is also the reason why in the proof of Theorem 14 we do not consider 
the more simply built compact linearly ordered spaces
\sp
\cl{$\;\tilde Y[S]\;\,:=\;\,\bigcup\limits_{\xi\in S}
\{\xi\}\!\times\![0,\omega^\xi]\!\times\! K 
\;\cup\;\{(\kappa,0)\}\!\times\!K$}
\sp
because for $\,K=[0,1]\,$ and $\,\kappa=\aleph_0\,$
all spaces $\;\tilde Y[S]\;(S\in{\cal S}_\kappa)\;$ are homeomorphic.
(If for infinite $\,S\subset\N\,$
every component of $\,\tilde Y[S]\,$ is condensed to 
a singleton then we obtain a well-ordered set 
order-isomorphic to $\,[0,\omega^\omega]\,$
since $\;\big(\sum_{\xi\in S}(\omega^\xi+1)\big)+1\,=\,\omega^\omega+1\,$.)
\mp\sp
In the proof of Theorem 14 the
the spaces $\,Y[S]\,$ are not necessarily incomparable. 
For example, if $\,K=[0,1]\,$ and $\,\kappa=\aleph_0\,$ 
then for arbitrary sets $\,S_1,S_2\in{\cal S}_\kappa\,$ the space
$\,Y[S_1]\,$ is obviously embeddable into $\,[0,1]\,$ and hence
into $\,Y[S_2]\,$. Fortunately this causes no problem 
in the following proof of Theorem 8.
\mp\sp
We verify Theorem 8 by expanding the previous proof 
of Theorem 14 with its notations.
Assume $\,\lambda\leq \kappa\,$ and $\,\kappa^{\aleph_0}=\kappa\,$
and put $\,\mu=\nu(\lambda,\kappa)\,$.
Let $\,K\,$ be a balanced linearly ordered continuum with $\,w(K)=\kappa\,$
and $\,|K|=\lambda^\mu\,$. (Such a space $\,K\,$ exists due to 
the basic constructions in Section 4.)
Consider the family 
$\;{\cal Y}[K,\kappa]\,=\,\{\,Y[S]\;|\;S\in{\cal S}_\kappa\,\}\;$
as defined above.
For $\,S\in{\cal S}_\kappa\,$ 
in the family of $\,{\cal C}[S]\,$ all components of
$\,Y[S]\,$  we select a subfamily as follows.
Let $\,{\cal C}[S]^*\,$ be the family of all 
$\,C\in{\cal C}[S]\,$ such that there is a component $\,D\in{\cal C}[S]\,$
with the property that $\;\max C\prec \min D\;$
and $\;\{\,x\in Y[S]\;|\;\max C\prec x\prec\min D\,\}=\emptyset\,$.
(Here $\,\preceq\,$ denotes the linear ordering of $\,Y[S]\,$.)
Of course, for every $\,C\in{\cal C}[S]^*\,$
there is precisely one component $\,D\in{\cal C}[S]\,$
with this property 
and we write $\,\rho(C)\,$ for this component $\,D\,$.
Obviuously, $\;|{\cal C}[S]^*|=|{\cal C}[S]|=\kappa\,$.
\sp
If $\,\Lambda\,$ is a linearly ordered set then 
$\,(a,b)\,$ is a {\it pair of consecutive points}
if and only if $\,b\,$ is the immediate successor of $\,a\,$
or, equivalently, if $\,a\,$ is the immediate predecessor of $\,b\,$.
So if $\,C,D\in{\cal C}[S]\,$ then 
$\,(C,D)\,$ is a pair of consecutive points in the canonically linearly
ordered family $\,{\cal C}[S]\,$ if and only if 
$\,C\in{\cal C}[S]^*\,$ and $\,D=\rho(C)\,$.
Clearly, if $\,C\in{\cal C}[S]\,$
has no immediate predecessor then
either $\,\min C=\min Y[S]\,$ (whence $\,C\,$ has no predecessors)
or $\,\min C\not=\min Y[S]\,$ and 
$\;\min C\,=\,\sup\,\{\,x\in Y[S]\;|\;x\prec \min C\,\}\,$.
And if a component $\,C\in{\cal C}[S]\,$
has no immediate successor then
either $\,\max C=\max Y[S]\,$ (whence $\,C\,$ has no successors)
or $\,\max C\not=\max Y[S]\,$ and 
$\;\max C\,=\,\inf\,\{\,x\in Y[S]\;|\;x\succ \max C\,\}\,$.
\sp
Now let 
$\;{\cal Q}(S)\,=\,\{\,(C,\rho(C))\;|\;C\in {\cal C}[S]^*\,\}\,$.
Let $\,H\,$ be a balanced 
linearly ordered continuum with $\,w(H)=|H|=\kappa\,$.
(Such a space $\,H\,$ exists because, referring to Section 4,
if $\,L\,$ is the long line of size and weight $\,\kappa\,$ then 
$\,H:=L^{[0,\omega[}\,$ fits because $\,\kappa^{\aleph_0}=\kappa\,$.)
Consider $\,H\,$ equipped with a linear ordering 
which generates the topology of the continuum $\,H\,$.
From the linearely ordered set $\,Y[S]\,$
we create a linearely ordered set $\,X[S])\,$
by filling the gap between $\,C\,$ and $\,D\,$ 
with an order-isomorphic copy of 
$\;H\setminus\{\min H,\max H\}\;$
for every pair $\,(C,D)\in {\cal Q}(Y[S])\,$.
Formally we put 
\sp
\cl{$X[S]\;=\;Y[S]\times\{\min H\}\,\cup\,
\bigcup\limits_{C\in{\cal C}[S]^*}\{\max C\}\times (H\setminus\{\max H\})$}
\sp
and consider $\,X[S]\,$ equipped with the lexicographic 
ordering. In view of Lemma 1 it is clear that 
the order topology of $\,X[S]\,$ is compact and connected.
So $\,X[S]\,$ is a linearly ordered continuum.
A moment's reflection suffices to see that 
$\,Y[S]\times\{\min H\}\,$ is a compact subspace 
of $\,X[S]\,$ homeomorphic with $\,Y[S]\,$
and that the lexicographic 
ordering of $\,X[S]\,$ restricted to
$\,Y[S]\times\{\min H\}\,$ generates the compact subspace topology. 
\sp
We have $\;|X[S]|=\lambda^\mu\;$ since 
$\;\kappa=|{\cal C}[S]|=|H|\leq |K|=\lambda^\mu\,$.
Furthermore, $\,w(X[S])=\kappa\,$ since $\,w(K)=w(H)=\kappa\,$
and $\,|{\cal C}[S]|=\kappa\,$.
Moreover, $\,\hat w(X[S])=\kappa\,$ since 
$\,\hat w(K)=\hat w(H)=\kappa\,$. 
The building blocks $\,K\,$ and $\,H\,$ 
of $\,X[S]\,$ are balanced of equal weight
but of different size. More precisely, 
the size of $\,H\,$ is smaller than the size of $\,K\,$.
This has the essential consequence that 
no interval in $\,K\,$ can be embedded into $\,H\,$
and no interval in $\,H\,$ can be embedded into $\,K\,$.
\mp
In order to conclude the proof of Theorem 8
let $\;S_1,S_2\in{\cal S}_\kappa\;$ 
and assume that $\,f\,$ is a homeomorphism from $\,X[S_1]\,$ 
onto a subspace of $\,X[S_2]\,$.
Of course, $\,f(X[S_1])\,$
is an interval in $\,X(S_2,H_\mu)\,$.
Since no interval in $\,K\,$ can be embedded into $\,H\,$
and no interval in $\,H\,$ can be embedded into $\,K\,$,
by similar arguments as in the proof of Theorem 4
we conclude that $\,f\,$ induces   
a monotonic function $\,\tilde f\,$ from the linearly ordered set $\,Y[S_1]\,$
into the linearly ordered set $\,Y[S_2]\,$ 
such that $\,\tilde f(C)\in {\cal C}[S_2]\,$
for every $\,C\in{\cal C}[S_1]\,$.
\sp
If $\,\xi\in S_1\,$ then $\,\tilde f\,$ 
maps the labeled component $\,\{z_\xi\}\times K\subset Z_\xi(K)\,$  
onto the labeled component 
$\,\{z_{\xi'}\}\times K\subset Z_{\xi'}(K)\,$  
for some $\,\xi'\in S_2\,$.
Moreover, if $\,\xi\in S_1\,$ and $\,\xi\not=\min S_1\,$
then it is inevitable that for every $\,\xi\in S_1\,$ 
the monotonic function $\,\tilde f\,$ maps the 
the whole linearly ordered block $\,Z_\xi(K)\,$ onto the 
whole linearly ordered block $\,Z_{\xi'}(K)\,$
for some $\,\xi'\in S_2\,$, whence $\,\xi'=\xi\,$.
(This would be not necessarily be true
if $\,\tilde f\,$ were only an embedding from 
the space $\,Y[S_1]\,$ into the space $\,Y[S_2]\,$!)
If $\,\xi=\min S_1\,$ then at least we can be sure that
$\,\tilde f(Z_\xi(K))\,$ is an order-isomorphic copy 
of $\,Z_{\xi}(K)\,$ within $\,Z_{\xi'}(K)\,$
such that $\,\tilde f(\{z_\xi\}\times K)=\{z_{\xi'}\}\times K\,$
for some $\,\xi'\,$ and this implies $\,\xi'=\xi\,$ as well.
(The case $\,\xi=\max S_1\,$ does not occur
because $\,S_1\,$ has no maximum due to
$\,|S_1|=\kappa\,$ and $\,S_1\subset[0,\kappa[\,$.)
So we conclude that $\,S_1\subset S_2\,$. 
This finishes the proof because by virtue of the following basic lemma
we can select a subfamily $\,{\cal S}_\kappa^*\,$ of $\,{\cal S}_\kappa\,$
such that $\,|{\cal S}_\kappa^*|=2^\kappa\,$ and 
$\,S_1\not\subset S_2\,$ whenever $\,S_1,S_2\in{\cal S}_\kappa^*\,$
and $\,S_1\not=S_2\,$. 
\mp
{\bf Lemma 12.} {\it If $\,Y\,$ is an infinite set of size $\,\kappa\,$
then there exists a family $\,{\cal Y}\,$ of subsets 
of $\,Y\,$ such that $\,|{\cal Y}|=2^\kappa\,$ and 
$\;|A\setminus B|=\kappa\;$ whenever 
$\;A,B\in{\cal Y}\;$ and $\,A\not=B\,$.} 
\mp  
{\it Proof.} Let $\,X\,$ be a set with $\,\kappa=|X|\,$.
Write $\,X=R\cup S\,$
with $\,|R|=|S|=|X|\,$ and $\,R\cap S=\emptyset\,$
and let $\,f\,$ be a bijection from $\,R\,$ onto $\,S\,$.
Put
$\;{\cal X}\,:=\,\{\,T\cup(S\setminus f(T))\;|\;T\subset R\,\}\,$.
Obviously, $\,|{\cal X}|=2^\kappa\,$ and 
$\;A\not\subset B\;$ whenever 
$\;A,B\in{\cal X}\;$ and $\,A\not=B\,$.
Consequently, the set $\;(A\times X)\setminus(B\times X)\,=\,
(A\setminus B)\times X\;$ is nonempty and hence of size $\,\kappa\,$
whenever 
$\;A,B\in{\cal X}\;$ and $\,A\not=B\,$.
Therefore, if $\,g\,$ is a bijection from $\;X\times X\;$ 
onto $\,Y\,$ then  
$\;{\cal Y}\,:=\,\{\,g(S\times X)\;|\;S\in{\cal X}\,\}\;$
is a family that fits, {\it q.e.d.}
\bp
{\it Remark.} Since $\,|I|<|J|\,$ for every interval $\,I\,$ in $\,H\,$ 
and every interval $\,J\,$ in $\,|K|\,$ it is straightforward 
to recover the compact linearly ordered
space $\,Y[S]\,$ from the linearly ordered continuum $\,X[S]\,$
for every $\,S\in{\cal S}_\kappa\,$.
Therefore, the statement about the family $\,{\cal G}_\lambda\,$
in Theorem 3 can be derived from Theorem 14 {\it provided that 
$\,\kappa^{\aleph_0}=\kappa\,$.} However, $\,\kappa^{\aleph_0}=\kappa\,$
is not assumed in Theorem 3 and the main part of the proof
of Theorem 3 is the construction of the family $\,{\cal F}_\kappa\,$.
\mp
It is worth mentioning that the previous proof of Theorem 8
shows a bit more. Obviously, we did not really use the fact
that $\,|K|=\lambda^{\nu(\lambda,\kappa)}\,$.
What we used is that $\,K\,$ is a balanced linearly ordered 
continuum of weight $\,\kappa\,$ and size {\it greater than} $\,|H|\,$.
In particular, the assumption $\,\kappa^{\aleph_0}=\kappa\,$
is superfluous as long as for the balanced building block $\,H\,$
we have $\,w(H)=\kappa\,$ and 
$\,|H|\leq \kappa^{\aleph_0}<|K|\,$. 
If we define $\,H:=L^{[0,\omega[}\,$
as above then $\,w(H)=\kappa\,$ and $\,|H|=\kappa^{\aleph_0}\,$.
Therefore, the proof of Theorem 8 also settles the following theorem. 
\mp
{\bf Theorem 15.} {\it If $\,\kappa^{\aleph_0}<\lambda\,$ 
and there exists a balanced linearly ordered continuum 
of weight $\,\kappa\,$ and size $\,\lambda\,$
then there exist $\,2^\kappa\,$ incomparable,  linearly 
ordered continua $\,X\,$ such that $\;\hat w(X)=w(X)=\kappa\;$
and $\,|X|=\lambda\,$.}
\bp
{\bff 10. On balanced, pathwise connected continua}
\mp
Referring to the notation of Section 2, if $\,X\,$
is a {\it balanced} linearly ordered continuum
then it is evident that the cone $\,\Psi(X)\,$ is 
a {\it balanced} continuum and $\,w(\Psi(X))=w(X)\,$
and  $\,|\Psi(X)|=|X|\,$. Therefore, by the same arguments as in the proof 
of Theorem 1, from Theorem 6 we immediately derive the following theorem.
\mp
{\bf Theorem 16.} {\it If 
$\;\max\{\mu^+,\lambda^{<\mu}\}\leq\kappa\leq \lambda^\mu\;$
then there exist $\,2^\kappa\,$ 
incomparable,  balanced, pathwise connected 
continua of weight $\,\kappa\,$
and size $\,\lambda^\mu\,$.}
\mp
As a consequence of (1.1) and Theorem 16 and Theorem 2, for every 
strong limit $\,\kappa\,$ 
there exist precisely $\,2^\kappa\,$ balanced, pathwise connected continua 
of weight $\,\kappa\,$ and size $\,2^\kappa\,$ up to homeomorphism. 
This enumeration result can be improved as follows.
\mp
{\bf Theorem 17.} {\it For every infinite cardinal $\,\kappa\,$ 
there exist precisely $\,2^\kappa\,$ pathwise connected, 
balanced continua 
of weight $\,\kappa\,$ and size $\,2^\kappa\,$ up to homeomorphism.}
\mp
The Hilbert cube $\,[0,1]^\kappa\,$
is a balanced pathwise connected continuum $\,X\,$ of weight $\,\kappa\,$ and 
size $\,2^\kappa\,$.
Therefore, Theorem 17 (and further enumeration theorems
about pathwise connected continua) 
can immediately be derived from (1.1) and the following theorem.
\mp\sp
{\bf Theorem 18.} {\it Let $\,J\,$ be a pathwise connected continuum
such that $\,w(J)\leq\kappa\,$.
Then there exist $\,2^\kappa\,$ mutually non-homeomorphic pathwise connected
continua $\,X_J\,$ of weight $\,\kappa\,$ and size $\,\max\{\kappa,|J|\}$.
All spaces $\,X_J\,$ are balanced if $\,J\,$ is balanced
and $\,w(J)=\kappa\,$.}
\mp
{\it Proof.} First of all we can assume that 
$\,J\setminus E\,$ is pathwise connected 
whenever $\,E\subset J\,$ and $\,|E|\leq 2\,$.
Because otherwise we replace 
$\,J\,$ with the product space $\,J\times J\,$.
(Clearly, $\,w(J\times J)=w(J)\,$ and $\,|J\times J|=|J|\,$
and if $\,J\,$ is balanced then $\,J\times J\,$ is balanced.) 
Furthermore assume, for the moment, that $\,\kappa>\aleph_0\,$.
By virtue of Theorem 3 there 
is a family $\,{\cal F}_\kappa\,$
of mutually non-homeomorphic {\it totally pathwise disconnected}
compact Hausdorff spaces of weight $\,\kappa\,$ and 
size $\,\max\{\kappa,\bc\}\,$
such that $\,|{\cal F}_\kappa|=2^\kappa\,$.
\mp
Fix any point $\,a\in J\,$ and 
consider the pathwise connected 
subspace $\;\tilde J\,:=\,J\setminus\{a\}\;$ of $\,J\,$.
For $\,Y\in {\cal F}_\kappa\,$ consider the compact product space
$\,Y\times J\,$. Similarly as in the proof of Theorem 3
let $\,X_J(Y)\,$ be the quotient 
of $\,Y\times J\,$ by its compact subspace $\,Y\times\{a\}\,$. 
So the subspace $\,Y\times \{a\}\,$ of $\,Y\times J\,$
is shrinked to a point $\,p\in X_J(Y)\,$ and 
the subspace $\;Y\times\tilde J\;$ of $\,Y\times J\,$
is identical with the subspace 
$\,X_J(Y)\setminus\{p\}\,$ of $\,X_J(Y)\,$.
Of course, $\,X_J(Y)\,$ is a compact Hausdorff space.
Since $\;|Y|=\max\{\bc,\kappa\}\;$ and $\,|J|\geq \bc\,$, we have 
$\;|X_J(Y)|=|Y\times J|=\max\{\kappa,|J|\}\,$. 
Since $\,Y\,$ is compact, by virtue of [2] 3.1.15
the filter of the neighborhoods of the compact subset 
$\,Y\times\{a\}\,$ in the space $\,Y\times J\,$ 
has a basis of size not greater than $\,w(J)\,$
and hence $\,p\,$ has a local basis of size not greater than $\,w(J)\,$
in the space $\,X_J(Y)\,$. Therefore, 
since $\,w(Y)=\kappa\,$ and $\,w(\tilde J)\leq w(J)\leq\kappa\,$,
we have $\;w(X_J(Y))=\kappa\,$. 
Of course, $\,X_J(Y)\,$ is pathwise connected.
Since every nonempty open subset of $\,X_J(Y)\,$
contains $\,\{y\}\times V\,$ for some $\,y\in Y\,$ 
and some open subset $\,V\,$ of $\,\tilde J\,$,
if $\,J\,$ is balanced and $\,w(J)=\kappa=w(Y)\,$
(and hence $\,|V|\geq |Y|\,$) then $\,X_J(Y)\,$ is balanced.
Since $\,Y\,$ is totally pathwise disconnected
and $\,\tilde J\setminus\{x\}\,$ is pathwise connected
for every $\,x\in\tilde J\,$,
the point $\,p\,$ is the unique point $\,x\,$ in $\,X_J(Y)\,$
such that $\,X_J(Y)\setminus\{x\}\,$ is not pathwise connected.
Hence the point $\,p\,$ is {\it topologically
labeled.} The path components of the space $\;X_J(Y)\setminus\{p\}\;$ 
are precisely the sets $\;\{y\}\times \tilde J\;$
with $\,y\in Y\,$. (For if $\,y_1,y_2\in Y\,$ 
and $\,z_1,z_2\in\tilde J\,$
and $\,f\,$ is a continuous injection from 
$\,[0,1]\,$ into $\;Y\times\tilde J\;$
with $\;f(0)=(y_1,z_1)\;$ and $\;f(1)=(y_2,z_2)\;$
then the projection $\,P\,$ of the pathwise connected continuum  $\;f([0,1])\;$
onto the basic space $\,Y\,$ is pathwise connected
and hence $\,|P|=1\,$ and hence $\,y_1=y_2\,$.) 
\mp
For distinct spaces $\,Y_1,Y_2\in{\cal F}_\kappa\,$
the spaces $\,X_J(Y_1)\,$ and $\,X_J(Y_2)\,$
are never homeomorphic because every $\,Y\in{\cal F}_\kappa\,$
can be recovered from $\,X_J(Y)\,$.
Indeed, define an equivalence relation on 
the topologically labeled subspace
$\,X_J(Y)\setminus\{p\}\,$ of $\,X_J(Y)\,$ 
such that two points are equivalent when they can be connected
by a path. Then the equivalence classes are the path components 
$\;\{y\}\times \tilde J\;$
of the space $\;X_J(Y)\setminus\{p\}\,=\,Y\times\tilde J\,$.
Consequently, the quotient space of $\,X_J(Y)\setminus\{p\}\,$ 
with respect to the equivalence relation is 
homeomorphic to $\,Y\,$. This concludes the 
proof of Theorem 18 for $\,\kappa>\aleph_0\,$.
\mp
Since any second countable pathwise connected continuum 
is balanced of size $\,\bc\,$ and $\,2^{\aleph_0}=\bc\,$,
the remaining case $\,\kappa=\aleph_0\,$ is settled 
by Theorem 2. However, it is much easier to track 
down non-homeomorphic spaces than incomparable spaces.
Moreover, the following easy construction shows that there 
exist $\,\bc\,$ mutually non-homeomorphic, 
second countable, {\it locally pathwise connected} continua.
(The clue in the proof of Theorem 2 is to destroy 
local pathwise connectedness at selected points.)
\sp
For distinct points $\,A,B\in\R^2\,$ let $\,[A,B]\,$ denote the closed straight
line segment which connects $\,A\,$ with $\,B\,$.
Let $\;{\cal M}\;$ be the family of all infinite sets of integers 
$\,k\geq 4\,$, whence  $\,|{\cal M}|=\bc\,$.
Now for each $\,M\in{\cal M}\,$ we define a set $\;A[M]\subset{\Bbb R}^2\;$ via
\sp
\cl{$ A[M]\;:=\;[0,1]\!\times\!\{0\}\;\cup\!
\bigcup\limits_{n\in M}\!\!\big([(2^{-n},0),(2^{-n},2^{-n})]\cup
\bigcup\limits_{k=1}^{n-1}
[(2^{-n},2^{-n}),((2^{-n}\!+\!3^{-nk},2^{-n+1})]\big)\,.$}
\sp
Obviously, $\,A[M]\,$ is a closed subset of the unit square $\,[0,1]^2\,$
and hence a compact subspace of $\,\R^2\,$.
Of course, the space $\,A[M]\,$ is pathwise connected.
Since the point $\,(0,0)\,$ causes no problems,
$\,A[M]\,$ is locally pathwise connected.
The $\,\bc\,$ spaces $\;A[M]\;(M\in{\cal M})\;$ are mutually non-homeomorphic 
because each set $\,M\in{\cal M}\,$ is completely determined 
by the topology of $\,A[M]\,$ via 
$\;M\,=\,\{\,p(x)\;|\;x\in A[M]\,\}\setminus\{1,2,3\}\;$
where $\,p(x)\,$ is the number of the path components 
of the space $\,A[M]\setminus\{x\}\,$.
\mp
{\it Remark.} The non-homeomorphic spaces $\,A[M]\,$ are far from being 
incomparable. Actually, it is straightforward to prove 
that if $\,M_1,M_2\in{\cal M}\,$ then $\,A[M_1]\,$ is 
homeomorphic to a subspace of $\,A[M_2]\,$. 
\bp
{\bff 11. On gaps between cardinals}
\mp
In this final section we present short proofs of the consistency results 
mentioned at the beginning of Section 1. 
For the sake of better reading put 
$\;\exp\kappa\,:=\,2^\kappa=\aleph_0^\kappa=\kappa^\kappa\;$ 
and $\,\nu(\kappa):=\nu(\kappa,\kappa)\,$
and $\,\rho(\kappa):=\nu(\aleph_0,\kappa)\,$,
whence $\;\nu(\kappa)\leq\rho(\kappa)\leq\kappa<\kappa^{\nu(\kappa)}
\leq \exp\rho(\kappa)\leq 2^\kappa\,$.
First we prove the following statement.
\mp\sp
{\bf Proposition 5.} {\it It is consistent with {\rm ZFC}
that the continuum function $\,\kappa\mapsto 2^\kappa\,$ is strictly increasing 
on the class of all cardinals and 

$\;\;${\rm (i)}
$\,|\{\,\lambda\;|\;\kappa<\lambda<2^\kappa\,\}|=2^\kappa\;$ 
for every regular $\,\kappa\,$
and 

$\;${\rm (ii)} $\,|\{\,\lambda\;|\;\kappa<\lambda<\kappa^{\nu(\kappa)}\,\}|=
\kappa^{\nu(\kappa)}\;$ whenever $\,\kappa=2^\mu\,$ for some $\,\mu\,$
and 

{\rm (iii)} 
$\,|\{\,\lambda\;|\;\kappa^{\nu(\kappa)}<\lambda<2^\kappa\,\}|=
2^\kappa\;$ whenever $\,\kappa=2^\mu\,$ for some singular $\,\mu\,$
with $\,\mu\not=\aleph_\mu\,$.}
\mp
{\it Proof.} In G\"odel's universe L, let $\,\Theta(\kappa)\,$ denote 
the smallest fixed point of the $\aleph$-function 
whose cofinality equals $\,\kappa^+\,$.
Then $\;|\{\,\lambda\;|\;\kappa<\lambda<\Theta(\kappa)\,\}|=\Theta(\kappa)\;$ 
holds in every generic extension of L.
By applying Easton's theorem [4, 15.18] one can 
create an Easton universe E generically extending L
such that the continuum function $\;\kappa\mapsto 2^\kappa=\kappa^+\;$
in L is changed into $\;\kappa\mapsto 2^{\kappa}=g(\kappa)\;$ in E
with $\,g(\kappa)=\Theta(\kappa)\,$ for every regular cardinal $\,\kappa\,$.
\sp
So in E we count $\,2^\kappa\,$ transfinite cardinals below $\,2^\kappa\,$
for every regular $\,\kappa\,$ which settles (i). Furthermore, 
since $\,\Theta(\cdot)\,$ is strictly 
increasing in L, in E we have 
$\,2^\alpha<2^\beta\,$ whenever $\,\alpha,\beta\,$ are {\it regular} cardinals
with $\,\alpha<\beta\,$. Therefore and 
since the {\it Singular Cardinal Hypothesis} 
\sp
\cl{$\;\exp{\rm cf}\,\kappa<\kappa\;\Longrightarrow\;
\kappa^{{\rm cf}\kappa}=\kappa^+\;$}
\sp
holds in every Easton universe
(see [1] Ex.~15.12) and by [1] Theorem 5.22, 
if $\,\mu\,$ is singular in E
then $\,2^\mu=(2^{<\mu})^+\,$ 
is a successor cardinal in E while 
if $\,\kappa\,$ is regular then 
$\,2^\kappa\,$ equals $\,\Theta(\kappa)\,$ which, of course, 
must be a limit cardinal. Consequently, in the universe E 
the continuum function 
is strictly increasing on {\it all} cardinals
and (i) is true. Thus in E we obtain (ii)
since $\,\kappa^\mu=\kappa=2^\mu<2^{(\mu^+)}=\kappa^{(\mu^+)}=\Theta(\mu^+)\,$ 
and (hence) $\,\nu(\kappa)=\mu^+\,$ in E.
Finally, since (i) and (ii) hold in E, in order to 
finish the proof by verifying (iii) in E, 
it is enough to make sure that $\,\mu\not=2^{<\mu}\,$
because $\,\kappa=2^\mu=(2^{<\mu})^+\,$ is regular and 
$\;2^{(\mu^+)}=\kappa^{\nu(\kappa)}\leq 2^\kappa=\exp((2^{<\mu})^+)\,$.
Well, in E the equation $\,\mu=2^{<\mu}\,$ is possible for a 
limit cardinal $\,\mu\,$ only if $\,\mu=\aleph_\mu\,$ because in E we have
$\;2^{<\mu}\,=\,\sup\{\,\Theta(\lambda^+)\;|\;\lambda<\mu\,\}\;$
and hence $\,2^{<\mu}\,$ is 
a fixed point of the $\aleph$-function, {\it q.e.d.}
\mp
Consider for a limit ordinal $\,\gamma\geq\omega_1\,$
the set 
$\;\Omega_\gamma\,:=\,\{\,\kappa\;\,|\,\;\aleph_0<\kappa<\aleph_{\gamma}
\,\land\;{\rm cf}\,\kappa=\aleph_0\,\}\;$ which is preserved 
as an (at least) uncountable set in all
generic extensions of G\"odel's universe L (cf.~[3] {\bf 13}).
By creating an Easton universe in which the 
function $\,\kappa\mapsto 2^\kappa\,$ is far from being 
strictly increasing, we can prove the following statement.
\mp
{\bf Proposition 6.} {\it It is consistent with {\rm ZFC}
that $\;|\{\,\lambda\;|\;\kappa^{\nu(\kappa)}<\lambda<\exp\rho(\kappa)\,\}|=
\exp\rho(\kappa)\;$
and 
$\;|\{\,\lambda\;|\;\exp\rho(\kappa)<\lambda<2^\kappa\,\}|=2^\kappa\;$
for every $\,\kappa\in \Omega_\gamma\,$.}
\mp
{\it Proof.} Choose ordinals $\,\alpha<\beta\,$ with
$\,\alpha=\aleph_{\alpha}\,$ and $\,\beta=\aleph_{\beta}\,$ 
and $\,{\rm cf}\,\alpha={\rm cf}\,\beta=\aleph_{\gamma+1}\,$.
Then 
$\;|\{\,\lambda\;|\;\alpha<\lambda<\beta\,\}|=\beta\;$
and $\;|\{\,\lambda\;|\;\kappa^+<\lambda<\alpha\,\}|=\alpha\;$
for every $\,\kappa\in \Omega_\gamma\,$.
Apply [3] 15.18 to create an Easton model of ZFC
where $\;\exp\aleph_0=\aleph_{1}\;$ and 
$\;\exp\aleph_{1}=\alpha\;$ and 
$\;\exp(\mu^{++})=\beta\;$
for every $\,\mu<\aleph_{\gamma}\,$. Then 
in our Easton model we have $\,2^\kappa=\beta\,$ and
$\,\exp\rho(\kappa)=\alpha\,$
for every $\,\kappa\in \Omega_\gamma\,$. 
Finally, if $\,\kappa\in \Omega_\gamma\,$ then 
in our Easton model we have $\,\kappa^{\nu(\kappa)}=\kappa^+\,$
due to the Singular Cardinal Hypothesis, {\it q.e.d.}
\bp\bp
{\bff References}
\mp
[1] J.~Baumgartner, {\it All $\aleph_1$-dense sets of reals can be 
isomorphic.} Fund.~Math.~{\bf 79} (1973) 101-106.
\sp
[2] R.~Engelking,  {\it General Topology, revised and completed edition.}
Heldermann 1989. 
\sp
[3] T.~Jech, {\it Set Theory}, 3rd ed. Springer 2002.
\sp
[4] I.~Juh\`asz, {\it Cardinal functions in topology}. 
Math.~Centrum Amsterdam 1983.
\sp
[5] G.~Kuba, {\it Counting metric spaces.}
Arch.~Math. {\bf 97} (2011) 569-578.
\sp
[6] G.~Kuba, {\it Counting linearly ordered spaces.}
Coll.~Math. {\bf 135} (2014) 1-14.
\sp
[7] G.~Kuba, {\it On compact subsets of the plane.} 
arXiv:2401.14985 [math.GN] (2024).
\sp
[8] K.~Kunen and J.E.~Vaughan, eds., {\it Handbook of set-theoretic topology.}
North Holland P.C.~1984.
\sp
[9] F.W.~Lozier and R.H.~Marty, {\it The number of continua.}
Proc.~American Math.Soc.~{\bf 40} (1973) 271-273.
\sp
[10] S.~Mazurkiewicz and W.~Sierpi\'nski, {\it Contribution … la topologie
des ensembles d'nom\-brables.} Fund.~Math.~{\bf 1} (1920) 17-27.
\sp
[11] J.G.~Rosenstein, {\it Linear orderings.} Academic Press 1982.
\sp
[12] L.A.~Steen and J.A.~Seebach, Jr., {\it Counterexamples in Topology.} 
Dover 1995.
\bp\bp
Gerald Kuba,  
Institute of Mathematics, University of Natural Resources and Life Sciences, 
1180 Wien, Austria. {\sl E-mail:} {\tt gerald.kuba(at)boku.ac.at}
\end